\documentclass{compositio}
\usepackage[latin1]{inputenc}
\usepackage[all]{xy}
\usepackage{amsfonts}
\usepackage{amsmath}
\usepackage{amssymb}
\usepackage{graphics}

\theoremstyle{plain}
\newtheorem{theo}{Theorem}[section]
\newtheorem*{theoA}{Theorem A}
\newtheorem*{theoB}{Theorem B}
\newtheorem{prop}[theo]{Proposition}

\newtheorem{lemma}[theo]{Lemma}

\theoremstyle{definition}
\newtheorem{defin}[theo]{Definition}
\newtheorem*{conj}{Conjecture}

\theoremstyle{remark}
\newtheorem{rque}[theo]{Remark}

\DeclareMathOperator{\A}{\mathbb{A}}
\DeclareMathOperator{\C}{\mathbb{C}}

\DeclareMathOperator{\F}{\mathbb{F}}
\DeclareMathOperator{\Q}{\mathbb{Q}}
\DeclareMathOperator{\R}{\mathbb{R}}
\DeclareMathOperator{\Z}{\mathbb{Z}}
\DeclareMathOperator{\n}{\mathbb{N}}

\DeclareMathOperator{\T}{\mathbb{T}}
\DeclareMathOperator{\bV}{\mathbb{V}}

\DeclareMathOperator{\cC}{\mathcal{C}}

\DeclareMathOperator{\cM}{\mathcal{M}}

\DeclareMathOperator{\cO}{\mathcal{O}}

\DeclareMathOperator{\cQ}{\mathcal{Q}}
\DeclareMathOperator{\cR}{\mathcal{R}}
\DeclareMathOperator{\cT}{\mathcal{T}}

\DeclareMathOperator{\cY}{\mathcal{Y}}
\DeclareMathOperator{\cX}{\mathcal{X}}

\DeclareMathOperator{\gc}{\mathfrak{c}}
\DeclareMathOperator{\gd}{\mathfrak{d}}
\DeclareMathOperator{\gm}{\mathfrak{m}}

\DeclareMathOperator{\gn}{\mathfrak{n}}
\DeclareMathOperator{\go}{\mathfrak{o}}

\DeclareMathOperator{\gq}{\mathfrak{q}}

\DeclareMathOperator{\gu}{\mathfrak{u}}

\DeclareMathOperator{\Ad}{Ad}
\DeclareMathOperator{\Ann}{Ann}

\DeclareMathOperator{\Cl}{\mathcal{C}}
\DeclareMathOperator{\End}{End}

\DeclareMathOperator{\Fitt}{Fitt}
\DeclareMathOperator{\Frob}{Frob}
\DeclareMathOperator{\GL}{GL}

\DeclareMathOperator{\G}{\mathcal{G}}
\DeclareMathOperator{\Hom}{Hom}
\DeclareMathOperator{\Ind}{Ind}
\DeclareMathOperator{\Nm}{\N_{F\!/\!\Q}}
\DeclareMathOperator{\N}{N}
\DeclareMathOperator{\PGL}{PGL}

\DeclareMathOperator{\SL}{SL}
\DeclareMathOperator{\SO}{SO}

\DeclareMathOperator{\Sym}{Sym}
\DeclareMathOperator{\Tam}{Tam}
\DeclareMathOperator{\Tor}{Tor}
\DeclareMathOperator{\tr}{tr}

\DeclareMathOperator{\ad}{ad}

\DeclareMathOperator{\bs}{\!\backslash\!}

\DeclareMathOperator{\full}{full}

\DeclareMathOperator{\id}{id}

\DeclareMathOperator{\limind}{\underset{\rightarrow}{lim}}

\DeclareMathOperator{\pr}{pr}
\DeclareMathOperator{\rf}{\mathrm{f}}
\DeclareMathOperator{\rE}{\mathrm{E}}
\DeclareMathOperator{\rH}{\mathrm{H}}
\DeclareMathOperator{\rL}{\mathrm{L}}

\newcommand{\cf}{{\it cf }}

\begin{document}

\title{On Ihara's lemma for Hilbert Modular Varieties}

\author{Mladen Dimitrov}
\email{dimitrov@math.jussieu.fr}
\address{Universit{\'e} Paris Diderot,  UFR de Math{\'e}matiques,  Case 7012, 2
place Jussieu,  75251 Paris, FRANCE}



\begin{abstract}
Let $\rho$ be a modulo $p$ representation of the absolute
Galois group of a totally real number field. Under the assumptions
that $\rho$ has  large image and admits a low weight crystalline 
modular deformation we show  that any low weight
crystalline deformation of $\rho$ unramified outside a
finite set of primes will be modular. We follow the  approach
of  Wiles  as generalized by Fujiwara. The main new ingredient is an Ihara type
lemma for the local component at $\rho$ of the middle
degree cohomology of a Hilbert modular variety. As an application we
relate the algebraic $p$-part of the  value at $1$ of the adjoint
$L$-function associated to a Hilbert modular newform to the
cardinality of the corresponding Selmer group.

\end{abstract}

\maketitle

\setcounter{tocdepth}{1}
\tableofcontents
\setcounter{tocdepth}{3}

\section{Introduction.}

\subsection{Statement of the main results.}

Let $F$ be a totally real number field of degree $d$,  ring of
integers  $\go$ and Galois closure $\widetilde{F}$.
Denote by  $J_F$  the set of   all embeddings of $F$ into  $\R$.
The  absolute Galois  group of a field $L$ is denoted  by $\G_L$.

\medskip
Let $f$ be a Hilbert modular newform over $F$ of level $\gn$ (an
ideal of $\go$), cohomological weight $k=\sum_{\tau\in J_F}k_\tau\tau$
($k_\tau\geq 2$ of the same parity) and put $k_0=\max\{k_\tau|\tau\in
J_F\}$.  For a prime $p$ 
and  an embedding $\iota_p: \overline{\Q}\hookrightarrow
\overline{\Q}_p$  one can associate
to $f$ and $\iota_p$ a  $p$-adic  representation (\cf \cite{Ta,Ta2}) :
\begin{equation}
\rho_{f,p}:\G_F \rightarrow \GL_2(\overline{\Q}_p),
\end{equation}
which is  irreducible, totally odd, unramified outside $\gn p$ and
characterized by the property that for each prime $v$ not dividing
$\gn p$ we have $\tr(\rho_{f,p}(\Frob_v))=\iota_p(c(f,v))$, where
$\Frob_v$ denotes a geometric Frobenius at $v$ and 
 $c(f,v)$ is the eigenvalue of $f$ for the standard Hecke operator $T_v$.
The embedding $\iota_p$ defines a  partition  $J_F=\coprod_v J_{F_v}$, 
where $v$ runs over the primes of $F$ dividing  $p$ and $J_{F_v}$
denotes the set  of embeddings of $F_v$ in $\overline{\Q}_p$. 
Then  $\rho_{f,p}|_{\G_{F_v}}$ is known to be de
Rham of   Hodge-Tate  weights  $(\frac{k_0-k_\tau}{2},
\frac{k_0+k_\tau}{2}-1)_{\tau\in J_{F_v}}$, 
unless  $k_0=2$, $\rho_{f,p}$ is residually 
reducible but not nearly-ordinary, $d$ is even and the automorphic
representation associated to $f$  is not a discrete series at any
finite place (\cf \cite{BR} and \cite{Ki}). 
If  $p>k_0$ is unramified in $F$ and relatively prime to $\gn$, then 
$\rho_{f,p}|_{\G_{F_v}}$ is crystalline (\cf \cite{breuil}).

\medskip
Such a $\rho_{f,p}$ is defined over the ring of integers $\cO$ of a
finite extension $E$ of $\Q_p$. Denote by 
 $\kappa$ be the residue field of $\cO$ and 
let $\overline{\rho}_{f,p}$  be the
semi-simplification of the reduction of  $\rho_{f,p}$ modulo a
uniformizer $\varpi$ of $\cO$.  We say that a two-dimensional
irreducible $p$-adic  (resp.  modulo $p$) representation of $\G_F$
is {\it modular} if it can be obtained by the above construction.
The following conjecture is a well known extension to an arbitrary
totally real field $F$ of a  conjecture of Fontaine and Mazur
\cite{FoMa} : 

\begin{conj}   A two-dimensional, irreducible, totally odd
  $p$-adic representation of $\G_F$ unramified
outside a finite set of primes  and  de Rham at all primes $v$ dividing
$p$ with distinct Hodge-Tate weights for each $F_v\hookrightarrow 
\overline{\Q}_p$,
  is modular, up to a twist by an integer power of the
 $p$-adic  cyclotomic character.
\end{conj}
Consider the assumption :

\medskip
\noindent $(\star)$ $p$ is  unramified in $F$ and  
$\textstyle p-1>\sum_{\tau\in J_F}(k_\tau-1)$.

\medskip
We provide some evidence for this conjecture by proving :
\begin{theoA}  Let  $\rho: \G_F \rightarrow  \GL_2(\overline{\F}_p)$ be a continuous representation.
Assume that :

\medskip
\noindent $\mathrm{{\bf (Mod}_{\rho}{\bf )}}$  there
exists  a Hilbert modular newform $f$ of  level prime to $p$
 and cohomological  weight $k$ as in $(\star)$, such
that $\overline{\rho}_{f,p}\cong \rho$, and

\medskip
\noindent   $\mathrm{{\bf (LI}_{Ind \rho}{\bf )}}$
the image of $\G_{\widetilde{F}}$  by 
$\otimes \Ind_F^{\Q}\rho=\!\!\underset{\tau\in \G_{\Q}/\G_F}{\bigotimes}\rho(\tau^{-1} \cdot \tau)$
is irreducible  of  order divisible by $p$.

\medskip
 Then all crystalline deformations of $\rho$  of weights
between $0$ and $p-2$  and unramified outside a finite  set of primes
are modular. 
\end{theoA}

\begin{rque} 
We have greatly benefited from the  work \cite{Fu} of Fujiwara, though
we use a  different approach (\cf \S\ref{strategy} for a more
detailed discussion). Furthermore, the proof of theorem A relies on
Fujiwara's results in the minimal case. 
Let us mention however that if  $P_{\rho}=\varnothing$  (\cf
Definition \ref{P})  
then  Theorem A is independent of the results of \cite{Fu}
(\cf Theorem \ref{minimal}).  
\end{rque}

\begin{rque} One can show that if $F$ is Galois over $\Q$ and if 
$f$ is a  Hilbert modular newform on $F$ which is not a theta series 
nor a twist of a base change of a Hilbert modular newform on $E\subsetneq F$, then
for  all but finitely many primes $p$, $\rho=\overline{\rho}_{f,p}$
satisfies  $\mathrm{{\bf (LI}_{Ind \rho}{\bf )}}$ for all 
$\iota_p: \overline{\Q}\hookrightarrow \overline{\Q}_p$.
\end{rque}

\begin{rque}\label{lowering}
The level lowering results of Jarvis \cite{Ja1, Ja2}, Fujiwara \cite{Fu2} and  
Rajaei \cite{Raj}, generalizing classical results of Ribet
\cite{ribet5} {\it et  al.} to the case of an arbitrary totally real
field $F$, imply that 
the newform $f$ in $\mathrm{{\bf (Mod}_{\rho}{\bf )}}$ can be chosen so that 
$\rho_{f,p}$ is a  minimally ramified deformation of $\rho$ in the sense of
Definition \ref{global}. 
\end{rque}

To a  Hilbert modular newform as above, Blasius and Rogawski \cite{BR}
attached a rank $3$ motive over $F$ with coefficients in
$\overline{\Q}$, pure of weight zero and autodual. For all $\iota_p$,
its  $p$-adic  realization $\Ad^0(\rho_{f,p})$ is given 
 by the adjoint action of $\G_F$ via $\rho_{f,p}$ on the space of two by two  trace zero
matrices. Denote by $\rL(\Ad^0(\rho_{f,p}),s)$ and $\Gamma(\Ad^0(\rho_{f,p}),s)$
the associated $\rL$-function and $\Gamma$-factor. 

In this setting, Beilinson and Deligne conjecture that 
the order of vanishing of $\rL(\Ad^0(\rho_{f,p}),s)$ at $s=1$ 
equals $\dim \rH^1_{\mathrm{f}}(F,\Ad^0(\rho_{f,p})\otimes \Q_p)-\dim
\rH^0(F,\Ad^0(\rho_{f,p})\otimes \Q_p)$, 
where  $\rH^1_{\mathrm{f}}$ is the Selmer
 group defined by Bloch and Kato (\cf \cite[\S2.1]{dfg}).
By a formula due to Shimura we know that $\rL(\Ad^0(\rho_{f,p}),1)$
is a non-zero multiple of the Peterson inner product of $f$, hence
does not vanish. Since $\rho_{f,p}$ is irreducible, by Schur's lemma
$\rH^0(F,\Ad^0(\rho_{f,p})\otimes \Q_p)=0$. Therefore, in our case, 
the Beilinson-Deligne conjecture is equivalent to  the vanishing of 
$\rH^1_{\mathrm{f}}(F,\Ad^0(\rho_{f,p})\otimes \Q_p)$. 

Let  $\Tam(\Ad^0(\rho_{f,p}))\subset \cO$ be the Tamagawa ideal
introduced   by Fontaine and Perrin-Riou (\cf \cite[I.4.1, II.5.3.3]{FoPR}). 

\begin{theoB}
Let $f$ be a Hilbert modular newform over $F$ of level prime to $p$ and
cohomological weight $k$ as in $(\star)$.  Let $\Omega_f^{\pm}\in
\C^\times\!\!/\cO^\times$ be any two 
complementary Matsushima-Shimura-Harder periods as in Definition
\ref{periods}.  If $\rho=\overline{\rho}_{f,p}$
satisfies $\mathrm{{\bf (LI}_{Ind \rho}{\bf )}}$, then 
\begin{enumerate}
\item $\displaystyle \left(\iota_p\left(\frac{\Gamma(\Ad^0(\rho_{f,p}),1)
\rL(\Ad^0(\rho_{f,p}),1)}{\Omega_f^+\Omega_f^-}\right)\right)_{\cO}=
\Tam(\Ad^0(\rho_{f,p}))
\Fitt_{\cO}\left(\rH^1_{\mathrm{f}}(F,\Ad^0(\rho_{f,p})\otimes
\Q_p/\Z_p)\right)$,

\item The Beilinson-Deligne conjecture holds  : 
$\rH^1_{\mathrm{f}}(F,\Ad^0(\rho_{f,p})\otimes \Q_p)=0$.
\end{enumerate}
\end{theoB}

An immediate corollary is that for $p$ as in the theorem, 
the $p$-adic valuation of $\Omega_f^+\Omega_f^-$ doesn't change when we twist $f$ 
by a Hecke character or when we choose different complementary
periods.

Theorem B is a first step towards the generalization to an
arbitrary totally real field  of the work \cite{dfg} of Diamond, Flach
and Guo  on the Tamagawa number  conjecture for $\Ad^0(\rho_{f,p})$
over $\Q$.
When $F$ is not $\Q$, it is an open problem how to identify 
the periods $\Omega_f^{\pm}$ used in Theorem B with the 
motivic periods attached to $f$ used in the formulation of
the Tamagawa number  conjecture.

\subsection{General strategy of the proof.}\label{strategy}

The method we use originates in  the work of
Wiles \cite{wiles2} and Taylor-Wiles \cite{TaWi}, later  developed by
Diamond \cite{diamond} and Fujiwara \cite{Fu}.

Let $\rho$ be as in Theorem A and  let $\Sigma$ be  finite  set of
primes  of $F$ not dividing $p$. In \S\ref{mini} we will define the
notion of a  $\Sigma$-ramified deformations   of $\rho$. 
By Mazur \cite{mazur2}  and  Ramakrishna \cite{Ram}, the functor
assigning  to a local complete noetherian $\cO$-algebra  $A$  with
residue field $\kappa$, the set of all  $\Sigma$-ramified
deformations   of $\rho$ to $A$, is representable by a
$\cO$-algebra $\cR_{\Sigma}$, called the  universal deformation ring.
Since $\rho$ is absolutely irreducible and odd, 
$\cR_{\Sigma}$ is topologically generated as a $\cO$-algebra
by traces of images of elements of $\G_F$ (\cf
\cite[pp.509-510]{wiles2}).  Moreover by the Cebotarev Density
Theorem, it is enough to take   traces of images of Frobenius elements
outside a finite set of primes.

Let  $S$ be a large finite set of primes and let  $\cT_{\Sigma}$ be
the $\cO$-subalgebra of 
$\prod_f \cO$ generated by $(\iota_p(c(f,v)))_{v\notin S}$ where $f$
runs over  all Hilbert modular newforms of weight $k$ such that
$\rho_{f,p}$ is a $\Sigma$-ramified deformation of
$\rho$. The $\cO$-algebra $\cT_{\Sigma}$ is local complete
noetherian and reduced. By the above discussion $\cT_{\Sigma}$ does not 
depend on the choice of $S$ and  the natural homomorphism
$\cR_{\Sigma}\rightarrow \prod_f \cO$  factors though a surjective
homomorphism of  local $\cO$-algebras
$\pi_{\Sigma}:\cR_{\Sigma}\rightarrow \cT_{\Sigma}$.   
Then Theorem A amounts to proving that $\pi_{\Sigma}$ is an
isomorphism. 

We follow Wiles'  method consisting in  showing  first  that
$\pi_{\varnothing}$ is an isomorphism  ({\it the minimal case}) and
then in proving, by induction on the cardinality of $\Sigma$, that
$\pi_{\Sigma}$  is an isomorphism  ({\it raising the level}). In
order to prove that $\cR_{\Sigma}$ is ``not too big'' we use
Galois cohomology via Proposition \ref{tangent}. In order to prove that
$\cT_{\Sigma}$ is ``not too small'' we  realize it {\it
geometrically}   as a local component of the  Hecke algebra acting
on the middle degree cohomology of some Shimura variety and then 
use this interpretation to study  congruences.

It is on that last point that our approach differs from Fujiwara's.
Whereas Fujiwara's uses some  quaternionic Shimura curves  or
Hida varieties  of dimension $0$,  we use the  $d$-dimensional  Hilbert
modular variety. The main ingredient in our approach is a result from
\cite{dim-hmv} guaranteeing the  torsion freeness of certain  local
components of the middle degree cohomology of a Hilbert modular
variety, which will be recalled in  the next section.

In the minimal case our modularity result is strictly included in Fujiwara's
since  we only treat the case $P_{\rho}=\varnothing$ (\cf Definition
\ref{P}) and furthermore we do not consider the ordinary
non-crystalline case. On the other hand our level raising results are
new, thanks to an Ihara type lemma for the middle degree cohomology of
Hilbert modular varieties (\cf Theorem \ref{ihara}). 
Our proof relies substantially on the $q$-expansion principle, which is
available  for Hilbert modular  varieties.  

Finally,  let us observe that whereas modularity lifting results
similar to Theorem A may be obtained in various ways (\cf
\cite{SkWi1,SkWi2,SkWi3},  \cite{Ta3} or \cite{Ki2}), the use of the
cohomology of Hilbert modular varieties seem to be inevitable to 
obtain results on the adjoint  $L$-functions and Selmer groups,
as  Theorem B.

\section{Cohomology of Hilbert modular varieties.}

In this section we state and prove a slightly more general 
version of a theorem in  \cite{dim-hmv}. We take
advantage of this opportunity to correct a wrong assumption in 
\cite{dim-hmv}, coming from a mistake in \cite{dimtildg}. We 
thank the referee for pointing out this error to us. 

\subsection{Hilbert modular varieties.} \label{hmv}

Denote by $\widehat{\Z}$ the profinite completion of $\Z$ and 
by $\A=(F\otimes \widehat{\Z}) \times (F \otimes_{\Q} \R)$ the ring of
adèles of $F$. For a prime $v$, let $\varpi_v$ denote an uniformizer
of $F_v$.

For an  open compact subgroup $U$ of $(\go\otimes\widehat{\Z})^\times$
we denote by $\Cl_{U}$ (resp. $\Cl_{U}^+$) the class group 
$\A^\times/F^\times U (F\otimes_{\Q}\R)^\times$
(resp. the narrow class group 
$\A^\times/F^\times U (F\otimes_{\Q}\R)_+^\times$, where
$(F\otimes_{\Q}\R)_+^\times$ denotes
the open cone of totally positive elements in
$(F\otimes_{\Q}\R)^\times$).

For an  open compact subgroup $K$ of $\GL_2(F\otimes\widehat{\Z})$
we denote by $Y_K$ the Hilbert modular  variety of level $K$
with complex points 
$\GL_2(F)\bs \GL_2(\A)/K\cdot\SO_2(F\otimes_{\Q}\R)(F\otimes_{\Q}\R)^\times$. 
By the Strong Approximation Theorem for $\GL_2$, the group of connected
components of $Y_K$ is isomorphic to $\Cl_{\det(K)}^+$. 

We will  consider the Hilbert modular  varieties as analytic varieties, 
except in the proofs of  Theorem \ref{ihara}
Proposition \ref{ihara-bis} and Proposition \ref{isom3} where we will
use integral models.

For an ideal    $\gn$ of $\go$, we consider the following open
compact subgroups of $\GL_2(F\otimes\widehat{\Z})$ : 

$$K_0(\gn)=\left\{\begin{pmatrix}a& b \\ c & d\end{pmatrix}\in
\GL_2(\go\otimes\widehat{\Z}) \, \Big{|} \, c\in \gn \right\} \text{ ,
} K_1(\gn)=\left\{\begin{pmatrix}a& b \\ c & d\end{pmatrix}\in 
K_0(\gn) \, \Big{|} \, d-1\in \gn \right\},$$

$$K_{11}(\gn)=\left\{\begin{pmatrix}a& b \\ c & d\end{pmatrix}\in
K_1(\gn) \, \Big{|} \, a-1\in \gn \right\}\text{ , and } 
K(\gn)=\left\{\begin{pmatrix}a& b \\ c & d\end{pmatrix}\in
K_{11}(\gn) \, \Big{|} \, b \in \gn \right\}.$$

For $?=0,1,11,\varnothing$ let $Y_?(\gn)$ be the Hilbert modular  variety of
level $K_?(\gn)$. 

Consider the following assumption :

{\bf (NT)  } $\gn$  does not divide $2$, nor $3$, nor $\Nm(\mathfrak{d})$.

In \cite[Lemma 1.4]{dimtildg} it is shown that under the assumption
{\bf  (NT)}, for all $x\in \GL_2(F\otimes\widehat{\Z})$, the group 
$\GL_2(F)\cap x K_1(\gn)x^{-1}(F\otimes_{\Q}\R)^\times\SL_2(F\otimes_{\Q}\R)$
is torsion free. This is not sufficient to claim that $Y_1(\gn)$ is
smooth. Here is  a   corrected statement : 

\begin{lemma}\label{torsionfree}
\begin{enumerate}
\item  $Y_K$ is smooth if, and only if,  for all $x\in
\GL_2(F\otimes\widehat{\Z})$, the quotient of the group 
$\GL_2(F)\cap x K x^{-1}(F\otimes_{\Q}\R)^\times\SL_2(F\otimes_{\Q}\R)$ by its center
is torsion free.

\item  If $\gn$ satisfies {\bf  (NT)}, then $Y_{11}(\gn)$ is smooth. 

\item  Let $\gu$  be a prime ideal of $F$ above a prime number $q$ such that :

$\bullet$ $q$ splits completely in $F(\sqrt{\epsilon} \,|\,
\epsilon\in\go^\times, \forall \tau\in J_F,  \tau(\epsilon)>0)$, and 

$\bullet$ $q\equiv -1 \pmod{4\ell}$ for all prime numbers $\ell$ such that
$[F(\zeta_\ell):F]=2$. 

Then  $Y_0(\gu)$ is smooth.

\item If $K'\subset K$ and $Y_K$ is smooth, then $Y_{K'}$ is
smooth and the natural morphism $Y_{K'}\rightarrow Y_K$ is étale with group 
$K/K'(K\cap \overline{F^\times})$. 
\end{enumerate}
\end{lemma}

\begin{proof}  The claims (i) and (iv) are well known, (ii) easily follows from \cite[Lemma 1.4]{dimtildg}. We will omit the proof of 
(iii) since it is very similar to the proof of lemma 
\ref{torsionfree-adjoint}(i) given below. 
\end{proof}

\medskip
From now on, we will only consider compact open subgroups $K$ factoring
as a product $\prod_v K_v$ over the primes $v$ of $F$, such that $K_v$ is
maximal for all primes $v$  dividing $p$ and   $Y_K$ is smooth. 
We denote by $\Sigma_K$ the set of primes $v$ where 
$K_v$ is not maximal.

\medskip
 For a $\cO$-algebra $A$, we denote by  $\bV_A$ the sheaf
of locally constant sections of 
\begin{equation}
\GL_2(F)\bs\left( \GL_2(\A)\times  V_{A}\right)/ K \cdot \SO_2(F\otimes_{\Q}\R)
(F\otimes_{\Q}\R)^\times \longrightarrow  Y_K,
\end{equation}
where  $V_A$ denotes the algebraic irreducible representation 
$\bigotimes_{\tau\in
  J_F}(\det^{\frac{k_0-k_\tau}{2}}\otimes\Sym^{k_\tau-2}A^2)$ of  
$\GL_2(A)^{J_F} \cong\GL_2(\go\otimes A) $ and 
$K$ acts on the right on $V_A$ via  its $p$-component 
$\prod_{v|p} K_v$.
Note that for $K'\subset K$, there is a natural projection 
$\pr: Y_{K'}\rightarrow Y_K$ and $\pr^*\bV_A=\bV_A$.  For
$g\in \GL_2(F\otimes\widehat{\Z})\cap M_2(\go\otimes\widehat{\Z})$ 
we define the Hecke correspondence $[KgK]$ on $Y_K$ by the usual diagram : 
\begin{equation}
\xymatrix@C=20pt@R=12pt{ & Y_{K\cap g^{-1}Kg} \ar[dl]_{\pr_1}\ar[r]^{\cdot g} & 
Y_{gKg^{-1}\cap K}\ar[dr]^{\pr_2} &  \\ Y_K & &  & Y_K}
\end{equation}

The  Hecke correspondences act naturally on the 
  Betti  cohomology groups $\rH^\bullet(Y_K,\bV_A)$ 
and on those with compact support $\rH_c^\bullet(Y_K,\bV_A)$.
If $K_v=\GL_2(\go_v)$,  we define the standard Hecke operators
$T_v=[K_v\left(\begin{smallmatrix}\varpi_v & 0 \\ 0 & 1 \end{smallmatrix}\right)K_v]$
and $S_v=[K_v\left(\begin{smallmatrix}\varpi_v & 0 \\ 0 &\varpi_v \end{smallmatrix}\right)K_v]$.

\subsection{Adjoint Hilbert modular varieties.} \label{adjoint-hmv}

For an  open compact subgroup $K$ of $\GL_2(F\otimes\widehat{\Z})$
we define the adjoint Hilbert modular  variety of level $K$ : 
\begin{equation}
Y_K^{\ad}=\GL_2(F)\bs \GL_2(\A)/\A^\times 
K\cdot\SO_2(F\otimes_{\Q}\R).
\end{equation}

Again, we have Betti  cohomology groups 
$\rH^\bullet(Y_K^{\ad},\bV_A)$ and Hecke action on them. 
In particular, if $K_v=\GL_2(\go_v)$, there is a Hecke operator
$T_v$ (the action of $S_v$ is trivial). 

We call $Y_K^{\ad}$ adjoint since it can be rewritten in terms of 
the adjoint group $\PGL_2$ as follows : 
\begin{equation}
Y_K^{\ad}=\PGL_2(F)\bs \PGL_2(\A)/
\overline{K}\cdot\mathrm{PSO}_2(F\otimes_{\Q}\R),
\end{equation}
where $\overline{K}$ is the image of $K$ in 
$\PGL_2(F\otimes\widehat{\Z})$.

The group of connected components of $Y_K$ is isomorphic to the 
quotient of  $\Cl_{\det(K)}^+$ by the image of $\A^{\times2}$, hence
it is a $2$ group. 
If $\det(K)=(\go\otimes\widehat{\Z})^\times$ then the group of connected
components of $Y_K$  is isomorphic to the
{\it narrow class group} $\Cl_F^+$ of $F$, while the group of connected
 of  $Y_K^{\ad}$ is isomorphic to the {\it genus group}
$\Cl_F^+/\Cl_F^2\cong \Cl_F^+/(\Cl_F^+)^2$. 
Each connected component of $Y_K^{\ad}$
can be defined more classically using the 
Hurwitz-Maass extension of the Hilbert modular group.

\begin{lemma}\label{torsionfree-adjoint}

\begin{enumerate}
\item  Let $\gu$  be a prime ideal of $F$ above a prime number $q$, such that  :

$\bullet$ $q$ splits completely in the ray class field of $F$ modulo $4$, and  

$\bullet$ $q\equiv -1 \pmod{4\ell}$ for all prime numbers $\ell$ such that $[F(\zeta_\ell):F)]=2$. 

Then  $Y_0^{\ad}(\gu)$ is smooth.
\item  If $K'\subset K$ and $Y_K^{\ad}$ is smooth, then $Y_{K'}^{\ad}$ is smooth and $Y_{K'}^{\ad}\rightarrow Y_K^{\ad}$ is étale with group 
$K/K'(K\cap {\A^\times})$. 
\end{enumerate}
\end{lemma}

\begin{proof}
We will show by contradiction that for all $x\in \GL_2(F\otimes\widehat{\Z})$, the quotient of the group $\GL_2(F)\cap x K_0(\gu) x^{-1}\A^\times \SL_2(F\otimes_{\Q}\R)$ by its center is torsion free. Suppose given an element $\gamma$ in that group which 
is torsion of prime order $\ell$ in the quotient. Consider the (quadratic) extension 
$F[\gamma]=F[X]/(X^2-\tr{\gamma}X+\det{\gamma})$ of $F$. Since
$\gamma_{\gu}\in K_0(\gu)F_{\gu}^\times$, it follows that $\gu$ splits in 
$F[\gamma]/F$.

If  $\ell$  is odd, then necessarily $F[\gamma]=F(\zeta_\ell)$. Our second assumption on 
$q$ implies then that $\gu$ is inert in $F[\gamma]$.  Contradiction. 

If  $\ell=2$, then $\tr{\gamma}=0$ and $\det{\gamma}\in F^\times\cap 
(\widehat{\Z}\otimes\go)^\times\A^{\times 2}$. By Class Field Theory,
the extension $F(\sqrt{\det{\gamma}})$ corresponds to a 
quotient of the class group $\Cl_{(1+4\widehat{\Z}\otimes\go)^\times}$, 
hence by our first assumption on $q$, $\gu$ splits in $F(\sqrt{\det{\gamma}})$.
On the other hand, by the second assumption $\gu$ is inert in 
$F(\sqrt{-1})$, hence $\gu$ is inert in $F(\sqrt{-\det{\gamma}})=F[\gamma]$.
Contradiction. 

This proves (i). The proof of (ii) is left to the reader. 
\end{proof}

\subsection{Freeness results.} 
Consider the maximal ideal $\gm_{\rho}=(\varpi,
T_v-\tr(\rho(\Frob_v)), S_v-\det(\rho(\Frob_v))\N_{F/\Q}(v)^{-1})$
of the abstract Hecke algebra  $\T^{S}=\cO[T_v, S_v |
\enspace v\notin S]$, where  $S$ is a finite set 
of primes containing $\Sigma_K\cup\{v\mid p\}$.
\begin{theo} \label{theo-AB} 
Let $K=\prod_v K_v \subset \GL_2(F\otimes\widehat{\Z})$ be an open compact
subgroup,  maximal for at primes $v$  dividing $p$ and 
such that  $Y_K$ is smooth.  Under  the assumptions  $(\star)$ and
  $\mathrm{\bf (LI_{\Ind\rho})}$ : 

\begin{enumerate}
\item 
$\rH_c^\bullet(Y_K,\bV_{\cO})_{\gm_{\rho}}=
\rH^\bullet(Y_K,\bV_{\cO})_{\gm_{\rho}}
=\rH^d(Y_K,\bV_{\cO})_{\gm_{\rho}}$ is a free $\cO$-module
of  finite rank.
\item 
$\rH^\bullet(Y_K,\bV_{E/\!\cO})_{\gm_{\rho}}
=\rH^d(Y_K,\bV_{E/\!\cO})_{\gm_{\rho}}$ is
 a  divisible $\cO$-module  of finite corank  and the Pontryagin pairing
$\rH^d(Y_K,\bV_{\cO})_{\gm_{\rho}} \times
\rH^d(Y_K,\bV_{E/\!\cO})_{\gm_{\rho}}\rightarrow E/\!\cO$ is a perfect
duality.
\end{enumerate}

Moreover, if $Y_K^{\ad}$ is smooth, then (i) and (ii) remain valid when we
replace $Y_K$ by $Y_K^{\ad}$. 
\end{theo}

\begin{proof} For $K=K_1(\gn)$ the theorem is proved in \cite[Theorems 4.4,
  6.6]{dim-hmv}, except that : 
  
$\relbar$  the assumption $\mathrm{\bf (LI_{\Ind\rho})}$ in \cite[\S 3.5]{dim-hmv}
is  formulated as follows : the restriction of $\rho$ to $\G_{\widetilde{F}}$   
is irreducible of  order divisible by $p$, and is not a twist by a 
character of any of its other $d-1$ internal conjugates. This is clearly 
implied by $\mathrm{\bf (LI_{\Ind\rho})}$. Conversely, if the assumption
 from \cite[\S 3.5]{dim-hmv} holds,
then by  \cite[Lemma 6.5]{dim-hmv} every irreducible 
$\G_{\widetilde{F}}$-representation annihilated by the characteristic 
polynomial of $(\otimes \Ind_F^{\Q}\rho)_{|\G_{\widetilde{F}}}$ is
isomorphic to $(\otimes \Ind_F^{\Q}\rho)_{|\G_{\widetilde{F}}}$,
 so in particular $(\otimes \Ind_F^{\Q}\rho)_{|\G_{\widetilde{F}}}$ is 
irreducible. Therefore these assumptions are equivalent.

$\relbar$  Theorem 4.4 is proved under the assumption {\bf (MW)}. However, 
this assumption is  only used through \cite[Lemma 4.2]{dim-hmv} and
under  the assumption  $\mathrm{\bf (LI_{\Ind\rho})}$ we  
can apply the stronger \cite[Lemma 6.5]{dim-hmv}, hence 
the results of \cite[Theorems 4.4]{dim-hmv} remain valid.

$\relbar$ we do not assume here  $\rho$ to be modular, 
but this doesn't affect the proofs of \cite{dim-hmv} in any way, since 
the arguments involving  $\rho$ are purely  group-theoretic.

\medskip
Let us now explain how these results extend to more general level
structures. Observe first that a  conjugate of $K$ has a normal subgroup of the
form $K(\gn)$ for some ideal $\gn\subset \go$. Hence a  conjugate of
$K$ contains $K_{11}(\gn)\cap K_{0}(\gn^2)$ as a normal subgroup. 
Therefore $Y_K$ admits a finite étale cover  isomorphic to $Y_{K_{11}(\gn)\cap
  K_{0}(\gn^2)}$,  and the latter has a finite abelian cover 
 $Y_{11}^1(\gn^2):=\coprod_{\gc}
M_1^1(\gc,\gn^2)$, where $\gc$ runs over a set of representatives
 of  $\Cl_{(1+\widehat{\Z}\otimes\gn^2)^\times}^{+}$ and 
$M_1^1(\gc,\gn^2)$ are the fine moduli spaces
defined in \cite[\S1.4]{dim-hmv}. The  following morphisms 
 of Hilbert modular varieties are  étale  :   
\begin{equation}\label{eq-exact-sequence}
\xymatrix@C=20pt@R=12pt{
 Y_{11}^1(\gn^2)\ar[r]  & Y_{11}(\gn^2)\ar[r]  & Y_{K_{11}(\gn)\cap
  K_{0}(\gn^2)}\ar[r]  &  Y_K \ar[r]  &  Y_K^{\ad}. }\end{equation}

Recall that each $M_1^1(\gc,\gn^2)$  is a fine moduli space admitting
an arithmetic model endowed with an universal Hilbert-Blumenthal
abelian variety. In \cite{dim-hmv,dimtildg} one proves 
various geometric results concerning $M_1^1(\gc,\gn)$,  such as 
the existence of minimal compactifications,  the existence of proper smooth toroidal
compactifications over $\Z_p$ and the extension of certain  vector bundles  
to these  compactifications, the construction of a
Berstein-Gelfand-Gelfand complex for distribution algebras over $\cO$,
having as consequence the  degeneracy at $\mathrm{E}_1$ of
the Hodge to De Rham  spectral sequence.  
By applying those constructions to each component of $Y_{11}^1(\gn^2)$, 
it follows that the highest  weight 
$\sum_{\tau\in J_F}(k_\tau-1)$ of $\otimes\Ind_F^{\Q}\rho$ does not occur in 
$\rH^i(Y_{11}^1(\gn^2)_{\overline{\Q}},\bV_{\kappa})$
for $i<d$. By \cite[Theorem 6.6]{dim-hmv} 
$\rH^i(Y_{11}^1(\gn^2)_{\overline{\Q}},\bV_{\kappa})_{\gm_{\rho}}$
vanishes for $i<d$ (it is important observe that the 
Hodge to De Rham  spectral sequence is $\T^S$-equivariant ; we refer
to  \cite[\S2.4]{dim-hmv} 
for a geometric definition of the Hecke correspondences). 

If $Y_{K''}\rightarrow Y_{K'}$ is an étale morphism of smooth Hilbert modular
varieties with group $\Delta$, the corresponding Hoschild-Serre
spectral sequence  is Hecke equivariant and yields
\begin{equation}\label{eq-HS}
\mathrm{E}_2^{j,i}=\rH^j(\Delta,
\rH^i(Y_{K''},\bV_{\kappa})_{\gm_{\rho}})\Rightarrow  
\rH^{i+j}(Y_{K'},\bV_{\kappa})_{\gm_{\rho}}.
\end{equation}
Starting from the vanishing 
of $\rH^i(Y_{11}^1(\gn^2),\bV_{\kappa})_{\gm_{\rho}}$ for $i<d$,
then applying (\ref{eq-HS}) to the morphisms of (\ref{eq-exact-sequence})
yields the vanishing of  $\rH^i(Y_K,\bV_{\kappa})_{\gm_{\rho}}$ 
and  $\rH^i(Y_K^{\ad},\bV_{\kappa})_{\gm_{\rho}}$) for $i<d$. The
theorem then follows by exactly the same 
arguments as in \cite[Theorems 4.4, 6.6]{dim-hmv}.
\end{proof}

\begin{prop}  \label{freeness} 
Suppose given an étale morphism  of smooth Hilbert modular
varieties $Y_K\rightarrow Y_{K'}$ with group
$\Delta$. Assume that $\Delta$ is an abelian $p$-group and that 
$\cO$ is large enough to contain the values of all its 
characters.  Then, under  the assumptions $(\star)$ and
  $\mathrm{\bf (LI_{\Ind\rho})}$,
$\rH^d(Y_K,\bV_{\cO})_{\gm_{\rho}}$ is a free $\cO[\Delta]$-module
and $\rH^d(Y_K,\bV_{\cO})_{\gm_{\rho}}
\otimes_{\cO[\Delta]}\cO \cong\rH^d(Y_{K'},\bV_{\cO})_{\gm_{\rho}}$
 as   $\T^{S}$-modules.
\end{prop}

\begin{proof}
By Theorem \ref{theo-AB}(i) $\rH^d(Y_K,\bV_{\cO})_{\gm_{\rho}}$
 is free over $\cO$, hence by Nakayama's lemma 
the desired freeness over $\cO[\Delta]$ is equivalent 
to the freeness of $\rH^d(Y_K,\bV_{\cO})_{\gm_{\rho}}\otimes\kappa$
over $\kappa[\Delta]$.  

Since $\Lambda:=\kappa[\Delta]$ is a local Artinial  ring, freeness is 
equivalent to flatness. 
Hence we have to show that 
$\Tor_i^{\Lambda}(\rH^d(Y_K,\bV_{\kappa})_{\gm_{\rho}}
,\kappa)=0$ for $i>0$ and
$\rH^d(Y_K,\bV_{\kappa})_{\gm_{\rho}}\otimes_{\Lambda}\kappa\cong  
\rH^{d}(Y_{K'},\bV_{\kappa})_{\gm_{\rho}}$.   

We reproduce here Fujiwara's {\it perfect complex} argument (\cf
\cite[Lemma 8.16]{Fu}) following the presentation of Mokrane and Tilouine
(\cf \cite[\S10]{MoTi}). 

\smallskip
Let $\cC^{\bullet}$ be the Godement resolution of the sheaf $\bV_{\kappa}$
on the (complex) variety $Y_K$. It has a natural action of
$\Lambda$ and there is an hypertor spectral sequence :
\begin{equation*}
\rE_2^{i,j}=\Tor_{-i}^{\Lambda}(\rH^j(\cC^{\bullet}),\kappa)
\Rightarrow \rH^{i+j}(\cC^{\bullet}\otimes_{\Lambda}\kappa).
\end{equation*}

By definition $\rH^j(\cC^{\bullet})=\rH^j(Y_K,\bV_{\kappa})$.
Since $Y_K\rightarrow Y_{K'}$ is étale  with group
$\Delta$,  it is a standard property of Godement's resolution that 
$\rH^j(\cC^{\bullet}\otimes_{\Lambda}\kappa)=\rH^j(Y_{K'},\bV_{\kappa})$ 
(\cf \cite[Lemma 8.18]{Fu}). Hence the spectral sequence becomes : 
\begin{equation*}
\rE_2^{i,j}=\Tor_{-i}^{\Lambda}(\rH^j(Y_K,\bV_{\kappa}),\kappa)
\Rightarrow \rH^{i+j}(Y_{K'},\bV_{\kappa}).
\end{equation*}
Since the Hecke operators are defined as correspondences, 
the spectral sequence is $\T^{S}$-equivariant and we can localize it
at $\gm_{\rho}$. By Theorem \ref{theo-AB}(i), we have
$\rH^d(Y_K,\bV_{\kappa})_{\gm_{\rho}}=0$, unless $j=d$. Therefore
the $\gm_{\rho}$-localization of 
the spectral sequence degenerates at $\rE_2$, and gives : 
\begin{equation*}\Tor_{-i}^{\Lambda}
(\rH^d(Y_K,\bV_{\kappa})_{\gm_{\rho}},\kappa)\cong  
\rH^{i+d}(Y_{K'},\bV_{\kappa})_{\gm_{\rho}}.
\end{equation*}
Another application of Theorem \ref{theo-AB}(i) yields 
$\rH^{i+d}(Y_{K'},\bV_{\kappa})_{\gm_{\rho}}=0$, unless $i=0$. 

Hence $\Tor_{-i}^{\Lambda}
(\rH^d(Y_K,\bV_{\kappa})_{\gm_{\rho}},\kappa)=0$, unless $i=0$
in which case 

$\rH^d(Y_K,\bV_{\kappa})_{\gm_{\rho}}\otimes_{\Lambda}\kappa=
  \Tor_{0}^{\Lambda}
(\rH^d(Y_K,\bV_{\kappa})_{\gm_{\rho}},\kappa)\cong
\rH^{d}(Y_{K'},\bV_{\kappa})_{\gm_{\rho}}$ as desired. 
\end{proof}

\subsection{Twisted Hilbert modular varieties and  Hecke
  operators.} \label{twisted}

Let $U$ be an  open compact subgroup of
$(\go\otimes\widehat{\Z})^\times$ and let  $K$ be an  open compact
subgroup of  $\GL_2(F\otimes\widehat{\Z})$ such that 
$K_{11}(\gn)\subset K \subset K_0(\gn)$, for some 
ideal $\gn\subset \go$.  Assuming that $U$ and $K$ decompose
as a product over all primes $v$, so does the group 
\begin{equation}
K'=\{x\in K | \det(x)\in U\}.\end{equation}

We define the twisted Hecke operators 
$T'_v=[K_v'\left(\begin{smallmatrix}\varpi_v & 0 \\ 0 & 1 \end{smallmatrix}\right)K_v']$
and $S'_v=[K_v'\left(\begin{smallmatrix}\varpi_v & 0 \\ 0 & \varpi_v
\end{smallmatrix}\right)K_v']$,  for  $v\nmid \gn$, and
 $U'_v=[K_v'\left(\begin{smallmatrix}\varpi_v & 0 \\ 0 & 1
\end{smallmatrix}\right)K_v']$, for  $v\mid \gn$. 

Note that if $v\notin \Sigma_{K'}$, then
$T'_v$, $S'_v$ and $U'_v$  coincide with the standard Hecke operators. 
In general, they depend on the choice of $\varpi_v$ in
the following way :  
if we replace  $\varpi_v$ by $\varpi'_v$ then $T'_v$ and $U'_v$
are multiplied by the invertible Hecke operator 
$U_{\delta}:=[K_v'\left(\begin{smallmatrix}\delta & 0 \\ 0 & 1
 \end{smallmatrix}\right)K_v']=K_v'\left(\begin{smallmatrix}\delta & 0 \\ 0 & 1
 \end{smallmatrix}\right)$,  with  $\delta=\frac{\varpi'_v}{\varpi_v}\in
 \go_v^\times$,  whereas $S'_v$ is  multiplied by its square.

For a finite order Hecke character $\psi$ of $\Cl_{K'\cap\A^\times}$, 
we denote by $[\psi]$ the $\psi$-isotypic 
part for  the   action of the Hecke operators $S_v$, $v\notin
\Sigma_{K'}$.

For a finite order character $\nu$ of
$(\go\otimes\widehat{\Z})^\times$, trivial on $U$,  we
denote by   $[\nu]$  the $\nu$-isotypic part for  the action of the Hecke 
operators $U_\delta$ for $\delta\in \go_v^\times$.

\subsection{Modified Poincaré pairing.}\label{poincare}

We keep  the notations from the previous paragraph we  denote by  
$\rH^d_!(Y_{K'},\bV_{\cO})$ the image of
$\rH^d_c(Y_{K'},\bV_{\cO})$ in  $\rH^d(Y_{K'},\bV_{\cO})$.
There is perfect  pairing $\langle\cdot,\cdot\rangle : V_{\cO}\times
V_{\cO}\rightarrow \cO$ described 
in \cite[p.270]{Ta}
such that for all $x\in (M_2(\cO)\cap \GL_2(E))^{J_F}$ we have 
$$\langle\cdot x,\cdot x\rangle=\det(x)^{(k_0-2)t}\langle\cdot,\cdot\rangle.$$
It  induces by cup product a pairing : 
$$\rH^d_!(Y_{K'},\bV_{\cO})\times \rH^d_!(Y_{K'},\bV_{\cO})\rightarrow  \rH^{2d}(Y_{K'},\bV_{\cO}).$$ 
By the above relation, for any character $\psi$ of
$\Cl_{K'\cap\A^\times}$, the pairing only relates 
the $ [\psi]$-part to the $[\psi^{-1}]$-part.  Moreover $\rH^d_!(Y_{K'},\bV_{\cO})[\psi]\cong 
 \rH^d_!(Y_{K'}^{\ad},\bV_{\cO}^\psi)$, where $\bV_{\cO}^\psi$ is the sheaf
locally constant sections of 
\begin{equation}
\GL_2(F)\bs \left(\GL_2(\A)\times  V_{\cO}^\psi\right)/\A^{\times(p)}
K'   \SO_2(F\otimes_{\Q}\R)\longrightarrow Y_{K'}^{\ad},
\end{equation}
 where the prime to $p$ idèles  $\A^{\times(p)}$ act on $V_{\cO}^\psi$
 via $\psi |\cdot |^{k_0-2}$   (note that this is compatible with the
 action of $K'$, since by definition $\psi$ is trivial on  $K'\cap
 \A^\times$). By fixing a character of the  connected components 
 of $Y_{K'}^{\ad}$ (note that this is a $2$-group and $p$ will be odd)  
we obtain a pairing 
\begin{equation}
[\enspace,\enspace] :   \rH^d_!(Y_{K'}^{\ad},\bV_{\cO}^\psi)\times 
\rH^d_!(Y_{K'}^{\ad},\bV_{\cO}^{\psi^{-1}})\rightarrow 
 \rH^{2d}(Y_{K'}^{\ad},\bV_{\cO}) \rightarrow \cO,\end{equation} 
which becomes perfect after extending scalars to $E$.

The adjoint of the Hecke operator $[K'xK']$ with respect to this pairing is 
$\psi(x) [K'\det(x)x^{-1}K']$. 
We will next replace this pairing by  a Hecke equivariant one.
Put $\iota=\left(\begin{smallmatrix} 0& -1\\ 
\gn & 0\end{smallmatrix}\right)$, where by an abuse of notation 
the ideal $\gn$ denotes also the corresponding idèle.   
Since we will be only interested in commutative Hecke algebras, 
we can restrict our attention to Hecke operators corresponding to diagonal
$x$ and in this case $\det(x)x^{-1}=\iota^{-1} x \iota$.

If $K$ is $K_{11}(\gn)$  or $K_{0}(\gn)$, then both $K$ and $K'$ are
normalized by $\iota$. If  $K_{11}(\gn)\subset K \subset K_{0}(\gn)$
is  of the  groups that will be considered in \S\ref{level}, one 
can easily check that for every diagonal $x$ one can choose the
$x_i$'s so that we simultaneously have 
$K' x K'=\coprod_i K' x_i$ and $K^\iota x K^\iota=\coprod_i K^\iota x_i$, 
where  $K_\iota'=\iota K' \iota^{-1}$. It follows that 
we have a commutative diagram :
\begin{equation}
\xymatrix@C=40pt@R=20pt{ \rH^d_!(Y_{K'}^{\ad},\bV_{\cO}^\psi)\ar[r]^{[K'\iota ]} 
\ar[d]^{[K'\iota^{-1} x \iota K']} & 
\rH^d_!(Y_{K_\iota'}^{\ad},\bV_{\cO}^\psi)\ar[r]^{\otimes \psi^{-1}} 
\ar[d]^{[K_\iota' x K_\iota' ]}&
\rH^d_!(Y_{K'}^{\ad},\bV_{\cO}^{\psi^{-1}}) \ar[d]^{\psi(x)[K'xK']} \\
 \rH^d_!(Y_{K'}^{\ad},\bV_{\cO}^\psi)\ar[r]^{[K'\iota ]}  & 
\rH^d_!(Y_{K_\iota'}^{\ad},\bV_{\cO}^\psi)\ar[r]^{\otimes \psi^{-1}} &
\rH^d_!(Y_{K'}^{\ad},\bV_{\cO}^{\psi^{-1}}) }. \end{equation}

By composing the above pairing with the first line in the diagram
we obtain a new pairing : 
\begin{equation}
 \langle\enspace,\enspace \rangle :  \rH^d_!(Y_{K'}^{\ad},\bV_{\cO}^\psi)\times 
\rH^d_!(Y_{K'}^{\ad},\bV_{\cO}^{\psi})\rightarrow  \cO,\end{equation} 
 that we call the {\it modified} Poincaré pairing.  It has 
  the advantage of being Hecke equivariant,  in particular $\T^S$-linear. 
 Under the assumptions of theorem \ref{theo-AB}(i) its $\gm_{\rho}$-localization yields
a  perfect $\cO$-duality of $\cO$-free modules : 
\begin{equation}
\langle\enspace,\enspace \rangle : 
\rH^d(Y_{K'}^{\ad},\bV_{\cO}^\psi)_{\gm_{\rho}}\times 
\rH^d(Y_{K'}^{\ad},\bV_{\cO}^{\psi})_{\gm_{\rho}}\rightarrow  \cO.\end{equation}

\section{Ihara's lemma for Hilbert modular varieties.}

Recall our  running  assumptions that  $K$  factors as a product $\prod_v
K_v$ over the primes $v$ of $F$, that $K_v$ is maximal for all primes
$v$  dividing $p$ and  that $Y_K$ is smooth.

Let $\gq$ be a prime not dividing $p$ and let  $S$ be a finite set 
of primes containing those dividing $p\gq$ and 
the set of primes $\Sigma_{K}$ where $K$ is not maximal.

Consider the maximal ideal $\gm_{\rho}=(\varpi,
T_v-\tr(\rho(\Frob_v)), S_v-\det(\rho(\Frob_v))\N_{F/\Q}(v)^{-1})$ 
of the abstract Hecke algebra  $\T^{S}=\cO[T_v,S_v | \enspace v\notin S]$. 
The Betti cohomology groups $\rH^d(Y_K,\bV_{\cO})$ defined in
\S\ref{hmv} are modules over $\T^S$. 

\subsection{Main theorem.}

Fix a finite index subgroup $U$ of $\go_{\gq}^\times$, and suppose that 
$K_{\gq}=\{x\in \GL_2(\go_{\gq}) | \det(x)\in U\}$.

In \S\ref{twisted} we defined Hecke operators 
$T'_{\gq}$, $S'_{\gq}$ (resp. $U'_{\gq}$, $U_{\delta}$, $\delta\in
\go_{\gq}^\times$) acting on 
$\rH^d(Y_K,\bV_A)$ (resp. on $\rH^d(Y_{K\cap K_0(\gq)},\bV_A)$).

Finally consider the  degeneracy maps 
 $\pr_1,\pr_2: Y_{K\cap K_0(\gq)} \rightarrow  Y_K$
used in the  definition of the Hecke correspondence $T'_{\gq}$.

\begin{theo} \label{ihara} 
Assume   that $(\star)$  and  
$\mathrm{\bf (LI_{\Ind \rho})}$ hold.
Then the $\gm_{\rho}$-localization of the 
 $\T^{S}$-linear  homomorphism :
$$\pr_1^*+\pr_2^*: \rH^d(Y_K,\bV_{\cO})^{\oplus 2}\rightarrow
\rH^d(Y_{K\cap K_0(\gq)},\bV_{\cO})$$
is injective with flat cokernel.
\end{theo}

\begin{proof}
Our proof is geometric and relies on the existence of smooth
models $\cY_K$ (resp. $\cY_{K\cap K_0(\gq)}$) of $Y_{K}$
(resp. $Y_{K\cap K_0(\gq)}$)  over an unramified extension of  $\Z_p$
and on the existence  of smooth toroidal compactifications
thereof. One should be careful to observe that $K\cap K_0(\gq)$ is  
maximal at primes dividing $p$. By the  Betti-étale comparison 
isomorphism the  cohomology groups 
$$W:=\rH^d(Y_{K,\overline{\Q}},\bV_{\kappa})_{\gm_{\rho}}
\text{ and } W_0(\gq): = 
\rH^d(Y_{K\cap K_0(\gq),\overline{\Q}},\bV_{\kappa})_{\gm_{\rho}},$$
are endowed with a structure of $\T^S[\G_{\Q}]$ modules. 
The  theorem is equivalent to the injectivity of
$\T^{S}[\G_{\Q}]$-linear homomorphism : 
 $$\pr_1^*+\pr_2^*: W^{\oplus 2}\rightarrow W_0(\gq).$$

The image of $\T^{S}_{\gm_{\rho}}$ in $\End_{\kappa}(W)$
is a local  Artinian ring and $(\gm_{\rho}^iW)_{i\geq 0}$ is a
finite decreasing filtration of $W$ by $\T^{S}[\G_{\Q}]$-modules. 
By  the torsion freeness result in Theorem \ref{theo-AB}(i), both
$W$ and the graded pieces $\gm_{\rho}^iW/\gm_{\rho}^{i+1}W$ are
quotients  of two $\T^{S}[\G_{\Q}]$-stable $\cO$-lattices in
$\rH^d(Y_{K,\overline{\Q}},\bV_{\cO})_{\gm_{\rho}}\otimes_{\cO} E$.
By a theorem of Brylinski  and Labesse  \cite{BL}, it follows
that the characteristic polynomial of $\otimes \Ind_F^{\Q}\rho$
annihilates the $\kappa[\G_{\Q}]$-module 
$\gm_{\rho}^iW/\gm_{\rho}^{i+1}W$ (\cf also \cite[Lemma 3]{diamond2}).  
It follows then from  $\mathrm{\bf (LI_{\Ind \rho})}$ and  \cite[Lemma
6.5]{dim-hmv}  
that every $\G_{\widetilde{F}}$-irreducible subquotient of $W$ 
is  isomorphic to $\otimes \Ind_F^{\Q}\rho$.  The same arguments 
apply also to $W_0(\gq)$.  Therefore we can check the above injectivity 
by checking it on the last graded pieces of the corresponding Fontaine-Laffaille
modules.

By Faltings' étale-crystalline comparison theorem and the degeneracy 
of the Hodge to De Rham spectral sequence (\cf  \cite[Theorem 5.13]{dim-hmv})
the claim would follow from the following lemma 
(although this part of the argument relies on the existence of toroidal 
compactifications of  $\cY_K$ and $\cY_{K\cap K_0(\gq)}$,  by  Köcher's Principle we 
can omit them as long as we are concerned with global sections
of the invertible bundle $\underline{\omega}^k\otimes\underline{\nu}^{1-k_0/2}$
(\cf  \cite[\S1.5,\S1.7]{dim-hmv})) : 
\end{proof}

\begin{lemma} The following  homomorphism is injective 
$$\pr_1^*+\pr_2^*:\rH^0(\cY_K{}_{/\kappa},\underline{\omega}^k\otimes\underline{\nu}^{1-k_0/2})^{\oplus
  2} \rightarrow
\rH^0(\cY_{K\cap K_0(\gq)}{}_{/\kappa},\underline{\omega}^k\otimes\underline{\nu}^{1-k_0/2}).$$
\end{lemma}

\begin{proof} Let $(g',g)$ be an element of the kernel : 
$\pr_1^*(g')=-\pr_2^*(g)$.
 
 Since  the homomorphism  is $U'_{\gq}$-equivariant for the
$U'_{\gq}$-action on the left hand side  given by the matrix
$\left(\begin{smallmatrix} T'_{\gq}  & -S'_{\gq}\Nm(\gq) \\ 1 &
    0\end{smallmatrix}\right)$, 
we may assume that $(g',g)$ is an eigenvector for $U'_{\gq}$. 
Similarly may assume that $g'$ is an eigenvector for
 $S'_{\gq}$. This implies that 
$g'$ is a multiple of $g$, hence  $\pr_2^*(g)=-\pr_1^*(g')$ is
a multiple of $\pr_1^*(g)$. On the other hand, $\pr_1^*(g)$ has
the same $q$-expansion as $g$, whereas the 
$q$-expansions of $\pr_2^*(g)$ and $g$ are related as follows
: for every  $x \in F\otimes \widehat{\Z}$,
\begin{equation}
c(\pr_2^*(g)|_{\gq},x)=\begin{cases} c(g,x\varpi_{\gq}^{-1}) & 
\text{ , if  } x_{\gq}\varpi_{\gq}^{-1}\in \go_{\gq} \text{, and  }\\
0 & \text{ , otherwise. }\end{cases}\end{equation}
It follows that  $c(g,x)=0$ for all $x$, which in vertu of 
the  $q$-expansion Principle implies  $g=0$. The proof of
Theorem \ref{ihara} is now complete.  
\end{proof}

\subsection{More cohomological results.}

Fix a finite index subgroup $U$ of $\go_{\gq}^\times$, and
suppose that $K_{\gq}= \{x\in K_1(\gq^{c-1})| \det(x)\in U\}$.

Consider the  degeneracy maps
\begin{equation}
\begin{split}
\pr_1,\pr_2:  Y_{K\cap K_1(\gq^{c})} \rightarrow
Y_{K\cap K_0(\gq^{c})} \rightarrow Y_{K} 
 \text{ and }\\
\pr_3,\pr_4: Y_{K\cap K_1(\gq^{c})\cap K_0(\gq^{c+1})}  
\rightarrow Y_{K\cap K_1(\gq^{c})}, 
\end{split}\end{equation}
used in the  definition of the Hecke correspondence $U'_{\gq}$ in
\S\ref{twisted}.

\begin{prop}  \label{ihara-bis}
Assume  that  $(\star)$  and  $\mathrm{\bf (LI_{\Ind \rho})}$ hold.
Then the $\gm_{\rho}$-localization of the 
 $\T^{S}$-linear  sequence :
$$ 0 \rightarrow
\rH^d(Y_{K\cap K_1(\gq^{c-1})},\bV_{\cO})
\overset{(\pr_1^*,-\pr_2^*)}{\longrightarrow}
\rH^d(Y_{K\cap K_1(\gq^{c})},\bV_{\cO})^{\oplus 2}
\overset{\pr_3^*+\pr_4^*}{\longrightarrow}
\rH^d(Y_{K\cap K_1(\gq^{c})\cap K_0(\gq^{c+1})},\bV_{\cO})$$
is exact and the last arrow has flat cokernel.
\end{prop}

\begin{proof}  We follow closely  Fujiwara's argument  
\cite[Proposition 5.13]{Fu2}, except for the last part of it where we
use a geometric argument instead (Fujiwara uses open compact 
subgroups which do not satisfy our running assumption to be maximal at
primes dividing $p$).    
 
It is enough to prove the exactness after tensoring with $\kappa$,
which by Theorem \ref{theo-AB}(i) amounts to  replacing $\bV_{\cO}$ by
$\bV_{\kappa}$. Put $K_0=K\cap K_1(\gq^{c-1})$,  $K_1=K\cap K_1(\gq^{c})$,

$$K_2=\begin{pmatrix} \varpi_{\gq} & 0 \\ 0 & 1 \end{pmatrix}(K\cap  K_1(\gq^{c}))
\begin{pmatrix} \varpi_{\gq}^{-1} & 0 \\ 0 & 1 \end{pmatrix}\text{ , and }$$

$$K_3= \begin{pmatrix} \varpi_{\gq} & 0 \\ 0 & 1 \end{pmatrix} 
(K\cap K_1(\gq^{c})\cap K_0(\gq^{c+1}))
\begin{pmatrix} \varpi_{\gq}^{-1} & 0 \\ 0 & 1 \end{pmatrix}=
K\cap K_1(\gq^{c})\cap  K^0(\gq) $$
where $K^0(\gq)=\begin{pmatrix} \varpi_{\gq} & 0 \\ 0 & 1 \end{pmatrix}
K_0(\gq)\begin{pmatrix} \varpi_{\gq}^{-1} & 0 \\ 0 & 1 \end{pmatrix}$ is
the opposite parahoric subgroup. 

For $i=0, 1,2,3$ put  $Y_i=Y_{K_i}$. By the above computations it is equivalent then to 
prove the exactness of the sequence :
$$ 0 \rightarrow
\rH^d(Y_0,\bV_{\kappa})_{\gm_{\rho}}
\overset{(\pr'_1{}^*,-\pr'_2{}^*)}{\longrightarrow}
\rH^d(Y_1,\bV_{\kappa})_{\gm_{\rho}}\oplus\rH^d(Y_2,\bV_{\kappa})_{\gm_{\rho}}
\overset{\pr'_3{}^*+\pr'_4{}^*}{\longrightarrow}
\rH^d(Y_3,\bV_{\kappa})_{\gm_{\rho}}\text{ , where }$$
$$\xymatrix@C=20pt@R=8pt{
& Y_3 \ar[dl]_{\pr'_3}\ar[dr]^{\pr'_4}\ar[dd]_{\pr}& \\ Y_1\ar[dr]_{\pr'_1}& &Y_2 \ar[dl]^{\pr'_2} \\& Y_0 & }$$
and the projections are induced by the inclusion of the open compact subgroups. 

\medskip

Taking  models of the $Y_i$ ($0\leq i\leq 3$) over $\Q$ and using
Betti-étale comparison isomorphisms turns the above sequence into a
sequence  of   $\T^{S}[\G_{\Q}]$-modules
$W_i:=\rH^d(Y_{i,\overline{\Q}},\bV_{\kappa})_{\gm_{\rho}}$. 
As in the proof of Theorem \ref{ihara}, the condition
$\mathrm{\bf (LI_{\Ind \rho})}$ implies that every
$\G_{\widetilde{F}}$-irreducible subquotient of $W_i$ ($0\leq i\leq 3$)
 is  isomorphic to $\otimes \Ind_F^{\Q}\rho$.
Therefore it is enough to check the exactness on the last graded pieces 
of the  Fontaine-Laffaille modules. This is the object of the
following : 
\end{proof}

\begin{lemma} The following sequence is exact : 
$$ 0 \rightarrow \rH^0({\cY_0}_{/\kappa},\underline{\omega}^k\otimes\underline{\nu}^{1-k_0/2})
\overset{(\pr'_1{}^*,-\pr'_2{}^*)}{\longrightarrow}
\rH^0({\cY_1}_{/\kappa},\underline{\omega}^k\otimes\underline{\nu}^{1-k_0/2})\oplus\rH^0({\cY_2}_{/\kappa},\underline{\omega}^k\otimes\underline{\nu}^{1-k_0/2})
\overset{\pr'_3{}^*+\pr'_4{}^*}{\longrightarrow}
\rH^0({\cY_3}_{/\kappa},\underline{\omega}^k\otimes\underline{\nu}^{1-k_0/2}).$$
\end{lemma}

\begin{proof}
We will adapt the analytic argument of \cite[Lemma 5.14]{Fu2} in order to 
show that  the coproduct $\cY_1\coprod_{\cY_3}\cY_2$ is  isomorphic to
$\cY_0$ as  $\kappa$-schemes.  

For $0\leq i\leq 3$, there exists  a fine
moduli scheme $\cY_i^1$ such that  $\cY_i^1\rightarrow \cY_i$ is a
finite \'etale   with group 
$$\Delta_i=\frac{F_+^\times\cap  \det(K_i)}
{(F^\times\cap  K_i)^2}, $$ where   
$\Delta_1=\Delta_2=\Delta_3\twoheadrightarrow \Delta_0$
(recall that by definition  $\cY_i^1$ has the same number of connected
components as $\cY_i$). Since $\cY_i^1\rightarrow \cY_0^1$ is 
$\Delta_i$-equivariant (where the action on  $\cY_0^1$ is via
the surjection $\Delta_i\twoheadrightarrow \Delta_0$), we have 
 $\cY_i\coprod_{\cY_i^1}\cY_0^1\cong \cY_0$  for all $i$. Hence  it is
 enough  to  show that  $\cY_1^1\coprod_{\cY_3^1} \cY_2^1\cong \cY_0^1$.  

We will show this claim using the following 
functorial description of the $\cY_i^1$'s :

\begin{enumerate}
\item $\cY_0^1$ classifies polarized Hilbert-Blumenthal abelian varieties
$A$  with  $\mu_{\gq^{c-1}}$-level structure $P$ and some additional
level structures that we will ignore since they are the same for
$\cY_i^1$ for all  $0\leq i\leq 3$ ;  
\item $\cY_1^1$ classifies polarized Hilbert-Blumenthal abelian varieties
$A$  with $\mu_{\gq^{c}}$-level structure $Q$ ; 
\item $\cY_2^1$ classifies polarized Hilbert-Blumenthal abelian varieties
$A$  with  a  $\mu_{\gq}$-subgroup $C$ and a 
 $\mu_{\gq^{c}}$-level structure $\overline{Q}$ in $A/C$ ;
\item $\cY_3^1$ classifies polarized Hilbert-Blumenthal abelian varieties
$A$  with $\mu_{\gq^{c}}$-level structure  $Q$ and a 
$\mu_{\gq}$-subgroup $C$ disjoint from the group generated by $Q$. 
\end{enumerate}

\noindent The morphisms $\pr'_j$ in the diagram above come from
forgetful functors described as follows :  
\begin{enumerate}
\item $\pr'_4(A,Q,C)=(A,Q\mod C, C)$, where 
$Q\mod C$ is  a $\mu_{\gq^{c}}$-level structure on $A/C$, since
$C$ is disjoint from the group generated by $Q$ ; 
\item $\pr'_3(A,Q,C)=(A,Q)$ ; 
\item  $\pr'_2(A,\overline{Q},C)=(A,\overline{Q}^{\gq})$ where 
it is important to note that $\overline{Q}^{\gq}$ is a well 
defined  $\mu_{\gq^{c-1}}$-level structure
 on  $A$  (not only in $A/C$) ; 
\item $\pr'_1(A,Q)=(A,Q^{\gq})$, where $Q^{\gq}$ 
is the $\mu_{\gq^{c-1}}$-level structure deduced from  $Q$ 
obtained by composing with the dual $\mu_{\gq^{c-1}}\hookrightarrow 
\mu_{\gq^{c}}$ of the natural projection $\go/\gq^{c}\rightarrow
\go/\gq^{c-1}$.  
\end{enumerate}

 We have  $\pr'_1\circ \pr'_3(A,Q,C)=(A,Q^{\gq})=(A,(Q\mod C)^{\gq})=
\pr'_2\circ \pr'_4(A,Q,C)$. 

We have to show that  given any two homomorphisms
$h_1 : \cY_1^1\rightarrow \cX$ and 
$h_2 : \cY_2^1\rightarrow \cX$ such that $h_1\circ \pr'_3=h_2\circ
\pr'_4$, there exists an unique homomorphism $h_0 : \cY_0^1\rightarrow \cX$
such that  $h_1=h_0\circ \pr'_1$ and  $h_2=h_0\circ \pr'_2$. 
By the  functorial description of the $\cY_0^i$'s and the $\pr'_j$'s
the claim is reduced to  
a simple lemma from group theory  saying that, if  $K_0$ is generated by
$K_1$ and $K_2$, then  the coproduct 
$K_0/K_1\coprod_{K_0/K_3} K_0/K_2=K_0/K_1\coprod_{K_0} K_0/K_2$ is a
singleton. 

Hence $\cY_1\coprod_{\cY_3}\cY_2\cong \cY_0$  yielding  an exact
sequence of  sheaves over    $\cY_0$ :   
$$ 0 \rightarrow\underline{\omega}^k\otimes\underline{\nu}^{1-k_0/2}\rightarrow
\pr'_1{}_*\pr'_1{}^*\underline{\omega}^k\otimes\underline{\nu}^{1-k_0/2}\oplus 
\pr'_2{}_*\pr'_2{}^*\underline{\omega}^k\otimes\underline{\nu}^{1-k_0/2}\rightarrow \pr_*\pr^*\underline{\omega}^k\otimes\underline{\nu}^{1-k_0/2}.$$
Since the functor of global sections is left-exact, this implies the lemma.
\end{proof}

\section{Twisting.}

Let $\rho: \G_F \rightarrow
\GL_2(\kappa)$ be a totally odd,  absolutely irreducible representation.

\subsection{Local twist types.}
For a prime $v$ of $F$,  we identify $\G_{F_v}$ with 
a decomposition subgroup
of $\G_F$ and denote by $I_v$ its inertia subgroup.
Let  $\rho_v$ be  the restriction of $\rho$ to $\G_{F_v}$. 
We normalize the  local Class Field Theory isomorphism so that 
the uniformizer $\varpi_v$ correspond to geometric Frobenius.

\smallskip
Over a totally real field $F$, twists of minimal conductor
exist locally, but not necessarily globally. This observation motivates 
the following definition, due to Fujiwara :

\begin{defin}\label{type}
Let $v$ be a prime of $F$ not dividing $p$.
A local twist type character for $\rho_v$ is a character
$\nu_v:\G_{F_v} \rightarrow \kappa^\times$ 
such that $\rho_v\otimes {\nu_v}^{-1}$ has minimal
conductor amongst  all twist of $\rho_v$ by characters of $\G_{F_v}$. 
For any prime $v$ we choose once for all a  local twist type character
$\nu_v$ and use the same notation for the character of $F_v^\times$
coming from local Class Field Theory. For simplicity, we 
choose $\varpi_v$ and $\nu_v$, so that $\nu_v(\varpi_v)=1$. 
\end{defin}

\begin{defin} \label{P}
Let $\Sigma_{\rho}$ be the set of primes $v$ not dividing $p$
such that $\rho_v\otimes {\nu_v}^{-1}$ is ramified. 

\smallskip
Let $S_{\rho}$
be the set of primes $v\in \Sigma_{\rho}$ such that $\rho_v$ is reducible. 

\smallskip
Let  $P_{\rho}$ be the set of primes $v\in \Sigma_{\rho}$ such that
$\rho_v$ is irreducible but $\rho_v|_{I_v}$ is reducible, and
$\Nm(v)\equiv -1 \pmod{p}$. 
\end{defin}

\medskip
Note that $\Sigma_{\rho}$, $S_{\rho}$ and $P_{\rho}$ do not change when
we twist $\rho$ by a character.

\subsection{Minimally ramified deformations.} \label{mini}

For a character $\mu$ taking  values in $\kappa^\times$,
we denote by $\widetilde{\mu}$ its Teichmüller lift.

Let $A$ be a local complete noetherian $\cO$-algebra with residue field
$\kappa$ and $\widetilde{\rho}_v:\G_{F_v} \rightarrow \GL_2(A)$  
be a lifting of $\rho_v$. For $F=\Q$, the following definition coincides 
 with the notion introduced in \cite{diamond3}.

\begin{defin}\label{mini-local}
We say that  $\widetilde{\rho}_v$  is a minimally ramified
if $\det\widetilde{\rho}_v|_{I_v}=\widetilde{\det\rho_v}|_{I_v}$ 
and additionally : 

$\relbar$ if  $v\notin\Sigma_{\rho}$, then 
$\widetilde{\rho}_v\otimes {\widetilde{\nu}_v}^{-1}$ is unramified. 

$\relbar$ if  $v\in S_{\rho}$, then 
$(\widetilde{\rho}_v\otimes \widetilde{\nu}_v^{-1})^{I_v}\neq 0$. 

$\relbar$ if  $v\in P_{\rho}$ and 
$(\rho_v\otimes \mu_v^{-1})^{I_v}\neq 0$  for some character 
$\mu_v:I_v\rightarrow \kappa^\times$, then 
$(\widetilde{\rho}_v\otimes \widetilde{\mu}_v^{-1})^{I_v}\neq 0$. 

\end{defin}

\begin{rque} \label{mini-rk}
\begin{enumerate}
\item  If  $\widetilde{\rho}_v$ is a  minimally ramified lifting of
$\rho_v$ then $\widetilde{\rho}_v \otimes \widetilde{\mu}$ is a
minimally ramified lifting of  $\rho_v\otimes
\mu$ for all characters $\mu:\G_{F_v}
\rightarrow \kappa^\times$.

\item  If $\widetilde{\rho}_v$ is a minimally ramified lifting of
$\rho_v$ then the Artin conductors of $\widetilde{\rho}_v$
 and $\rho_v$ coincide and $\det\widetilde{\rho}_v|_{I_v}$
is the Teichmüller lift  of $\det\rho_v|_{I_v}$.
The converse  holds if  $\rho_v$ has minimal
conductor among its twists and $v\notin P_{\rho}$  (\cf
 \cite[Remark 3.5]{diamond3}).
\end{enumerate}
\end{rque}

Let $\chi_p: \G_F \rightarrow \Z_p^\times$ be  the $p$-adic
cyclotomic character. 
Fix a finite $p$-power order character  $\phi: \G_F \rightarrow
\cO^{\times}$ of  conductor prime to $p$, and 
define $\psi: \G_F \rightarrow \cO^{\times}$ as the unique character
such that $\psi\phi^{-2}$  is the Teichmüller lift of
$(\chi_p\mod{p})\cdot\det\rho$.

\begin{defin}\label{global}
Let $\Sigma$ be  finite  set of primes  of $F$ not dividing $p$.
Let $A$ be a local complete noetherian  $\cO$-algebra  with residue
field $\kappa$. We say that a deformation $\widetilde{\rho}: \G_F
\rightarrow  \GL_2(A)$ 
of $\rho$ to $A$ is  $\Sigma$-ramified, if  the following three
conditions hold : 

$\relbar$  $\widetilde{\rho}\otimes \phi^{-1}$ is  minimally ramified
at all primes $v\notin \Sigma$, 
$v\nmid p$ (\cf Definition \ref{mini-local}),

$\relbar$ $\widetilde{\rho}$ is crystalline at each primes $v$
dividing $p$ with Hodge-Tate weights $(\frac{k_0-k_\tau}{2}, 
\frac{k_0+k_\tau}{2}-1)_{\tau\in J_{F_v}}$,

$\relbar$ $\det\widetilde{\rho}=\chi_p^{1-k_0}\psi $.

\smallskip
A $\varnothing$-ramified deformation is called minimally ramified.
\end{defin}
Note that if  $\rho_{f,p}$ is a $\Sigma$-ramified  deformation of
$\rho$, then the central character of  $f$ has to be
$\psi |\cdot|^{2-k_0}$. When $p$ is odd  every $p$-power
character of $\G_F$ has a square root, hence the determinant 
of any finitely ramified low weight crystalline
deformation of $\rho$ is of the form $\chi_p^{1-k_0}\psi $, 
for some $\psi$ as above.

\subsection{Auxiliary level structures.}

Under the assumption    $\mathrm{\bf (LI_{\Ind \rho})}$, a standard argument 
(\cf \cite[\S 12]{Ja2}) using the Cebotarev Density Theorem  implies that 
there exist infinitely many primes $\gu$ of  $F$ as in 
lemma \ref{torsionfree-adjoint}(i), such that
\begin{enumerate}
\item $\Nm(\gu)\not\equiv  1 \pmod{p}$ and 
\item  $\phi$ and $\rho$ are unramified at $\gu$, and $\tr(\rho(\Frob_{\gu}))^2
\not\equiv\psi(\gu)\Nm(\gu)^{k_0-2} (\Nm(\gu)+1)^2 \pmod{\varpi}$. 
\end{enumerate}

In particular this implies that $\rL_{\gu}(\Ad^0(\rho),1)\in\kappa^\times$. 
Let us fix such a prime $\gu$ and denote by $\alpha_{\gu}$ and  
$\beta_{\gu}$ the eigenvalues of $\rho(\Frob_{\gu})$. 

\begin{lemma}\label{aux}
The natural projection $\cT_{P_\rho\cup \{\gu\}}\rightarrow
\cT_{P_\rho}$ is an isomorphism. 
\end{lemma}

\begin{proof}
It amounts to proving that if $f$ is a newform of weight $k$, central 
character $\psi |\cdot|^{2-k_0}$  and level
prime to $p$, and if $\rho_{f,p}$ is a deformation of $\rho$ then the
local component $\pi_{\gu}$ of the associated automorphic representation $\pi$ 
is unramified. Since $\rho_{\gu}$ is unramified,
if  $\pi_{\gu}$ is ramified, then necessarily the valuation of its conductor is $1$ or $2$. 
Since $\pi_{\gu}$ has unramified central character this implies that 
$\dim \pi_{\gu}^{K_0(\gu)}=1$ or  $\dim \pi_{\gu}^{K_0(\gu^2)}=1$. In the first case $\pi_{\gu}$ is a  special representation, hence 
$\alpha_{\gu}\equiv \beta_{\gu}\Nm(\gu)^{\pm 1}\pmod{\varpi}$. 
In the second case $\pi_{\gu}$ is either a ramified principal series, 
in which case $\Nm(\gu)\equiv  1 \pmod{p}$, or a supercuspidal 
representation, in which case $\Nm(\gu)\equiv  -1 \pmod{p}$ and 
$\tr(\rho(\Frob_{\gu}))\equiv  0 \pmod{\varpi}$.
In both cases this contradicts our assumptions. 
\end{proof}

By lemmas \ref{torsionfree}(iii) and  \ref{torsionfree-adjoint}(i), for all
$K\subset K_0(\gu)$, $Y_K$ and $Y_K^{\ad}$ are smooth. 
However by lemma \ref{aux} the additional level at $\gu$ does not modify
the local components of the  Hecke algebras and cohomology modules
that we consider, hence we will omit it in our notations.

\subsection{Level structures and Hecke operators associated to
  $\rho$.} \label{level}  

The cohomology of the Hilbert modular varieties for the level structures 
that we will introduce in  this paragraph will play an important role
in the study of modular deformations of $\rho$. 

For  $v$ not dividing $p$ denote by   $c_v$ be the valuation of the Artin
conductor of $\rho_v\otimes{\nu_v}^{-1}$  and by 
 $d_v$ the dimension of $(\rho_v\otimes\nu^{-1})^{I_v}$(\cf
Definition \ref{type}). Put  $c_v=d_v=0$ if $v$ divides $p$. Define
\begin{equation}\begin{split}
K_v'=\ker(K_1(v^{c_v})\overset{\det}{\longrightarrow}\go_v^\times
\overset{\widetilde{\nu}_v\phi}{\longrightarrow} \cO^\times)\text{ ,
  and }\\
K_v''=\ker(K_1(v^{c_v})
\cap K_0(v^{c_v+d_v})\overset{\det}{\longrightarrow}\go_v^\times
\overset{\widetilde{\nu}_v\phi}{\longrightarrow} \cO^\times).
\end{split}
\end{equation} 
Note that for all but finitely many primes $v$, we have
$\nu_v|_{\go_v^\times}=\phi|_{\go_v^\times}=1$.   

For a prime $\gu$ as in lemma \ref{torsionfree-adjoint}(i)
and  a finite set of primes $\Sigma$ of $F$ not dividing $p$ we put 
 $\gn_{\Sigma}=\gu \prod_{v\in\Sigma} 
v^{c_v+d_v}\prod_{v\notin\Sigma} v^{c_v}$ and  
\begin{equation}
K_{\Sigma}=K_0(\gu)\cap\prod_{v\in\Sigma}K_v''\prod_{v\notin\Sigma} 
K_v'\subset K_0(\gn_{\Sigma}) \text{ ,  and } K_\rho=K_{\varnothing}. 
\end{equation}

As  in \S\ref{twisted} we define Hecke operators
$U_{\delta}:=K_v'\left(\begin{smallmatrix}\delta & 0 \\ 0 & 1
 \end{smallmatrix}\right)$ and 
$S_{\delta}:=K_v'\left(\begin{smallmatrix}\delta & 0 \\ 0 & \delta
 \end{smallmatrix}\right)$,  for all $v$ where
$\delta=\frac{\varpi'_v}{\varpi_v}\in  \go_v^\times$ ; 
$T'_v=[K_v'\left(\begin{smallmatrix}\varpi_v & 0 \\ 0
    & 1 \end{smallmatrix}\right)K_v']$ and $S'_v=[K_v'\left(\begin{smallmatrix}\varpi_v & 0 \\ 0 & \varpi_v
\end{smallmatrix}\right)K_v']$ for  $v\notin \Sigma$  such that 
$c_v=0$ ; $U'_v=[K_v'\left(\begin{smallmatrix}\varpi_v & 0 \\ 0
    & 1 \end{smallmatrix}\right)K_v']$ for  $v\notin \Sigma$ such that 
$c_v>0$ ;  $U''_v=[K_v''\left(\begin{smallmatrix}\varpi_v & 0 \\ 0
    & 1 \end{smallmatrix}\right)K_v'']$ for  $v\in \Sigma$.

\medskip
Let $Q$ be a finite set of primes as in \S\ref{R_Q} and let
$\Delta_{\gq}$ be the $p$-Sylow of $(\go/\gq)^\times$. Put     
\begin{equation}
K_{0,Q}=  K_\rho \cap \prod_{\gq\in Q} K_0(\gq) \text{ ,  and }
K^Q = K_\rho \cap \prod_{\gq\in Q}K^Q_{\gq}
\end{equation} 
where $K^Q_{\gq}$ is the kernel of the composition of 
 $ K_0(\gq)  \rightarrow (\mathfrak{o}/{\gq})^{\times},
 \left(\begin{smallmatrix} 
 a & b \\ c & d \end{smallmatrix}\right)\mapsto 
\frac{a_{\gq}}{d_{\gq}}$  with the natural projection
 $(\mathfrak{o}/{\gq})^{\times}
\rightarrow \Delta_{\gq}$. 

For   $\gq\in  Q$ and $\delta\in \gq$, the operator 
$S_{\delta}:=K^Q_{\gq}\left(\begin{smallmatrix}\delta & 0 \\ 0 & \delta
 \end{smallmatrix}\right)$ is trivial, the operator 
$U_{\delta}:=[K^Q_{\gq}\left(\begin{smallmatrix}\delta & 0 \\ 0 & \delta
 \end{smallmatrix}\right)K^Q_{\gq}]$ depends only on the 
image of $\delta$ in $\Delta_{\gq}$, and 
the operator $U_{\gq}:=[K^Q_{\gq}\left(\begin{smallmatrix}\varpi_{\gq} & 0
    \\ 0 & \delta   \end{smallmatrix}\right)K^Q_{\gq}]$
depends on the choice of $\varpi_{\gq}$ is described in \S\ref{twisted}.

\subsection{Decomposing the central action.}

Since our aim is to study automorphic forms with fixed 
central character, we will only consider open sugroups
$K\subset K_\rho$ such that $K\cap\A^\times=K_\rho\cap\A^\times$. 
Consider the idèle class group
\begin{equation}
\Cl_\rho:=\Cl_{K_\rho\cap\A^\times}. 
\end{equation}

The natural inclusions induce the  following commutative diagram, 
where all morphisms are étale for the indicated (abelian) groups : 
\begin{equation}
\xymatrix@C=80pt@R=20pt{Y_{\rho}\ar[d]_{\Cl_\rho}
\ar[r]^{F^\times K_1(\gn_{\varnothing})/F^\times K_\rho}
&Y_1(\gn_{\varnothing})\ar[d]^{\Cl_{1+\widehat{\Z}\otimes\gn_{\varnothing}}} \\
Y_{\rho}^{\ad}
\ar[r]^{\A^\times K_1(\gn_{\varnothing})/\A^\times K_\rho}&Y_{\gn_{\varnothing}}^{\ad} }
\end{equation}

If $v\in \Sigma_\rho$  the $p$-Sylow subgroup of $(\go/v)^\times$  injects naturally in $\A^\times K_1(\gn_{\varnothing})/\A^\times K_\rho$
({\it a fortiori} in $F^\times K_1(\gn_{\varnothing})/F^\times K_\rho$), hence
acts freely on $Y_{\rho}^{\ad}$ and $Y_{\rho}$. 
It follows that the étale morphism $Y_{\rho}\rightarrow Y_{\gn_{\varnothing}}^{\ad}$ 
factors through an étale morphism $Y_{\rho}\rightarrow Y_{\rho}^{\Delta}$
with group the $p$-group 
\begin{equation}
\Delta_\rho^\phi=(p\text{-Sylow of }\Cl_\rho) \times \prod_{v\in \Sigma_\rho} (\go_v^\times/\ker(\phi_v)).\end{equation}

Recall that   $[\psi]$ denotes the $\psi $-isotypic 
part for  the   action of the Hecke operators $S'_v$,
$v\notin\Sigma_\rho$, where $\psi $ is seen as a finite
order Hecke character of  $\Cl_\rho$, and 
that   $[\phi\widetilde{\nu}]$ denotes the intersection of the 
$\phi_v\widetilde{\nu}_v$-isotypic parts for  the action of the Hecke 
operators $U_\delta$ for $\delta\in \go_v^\times$.

For $v\notin \Sigma_\rho$ we have $U_\delta^2=S_\delta$ 
and since $p$ is odd, the $\phi_v$-action at those $v$ is 
determined by the  action of the central character.

Hence the $[\psi ,\widetilde{\nu}\phi]$-part is the intersection of the 
$[\phi^2,\phi]$-part for the action of the  $p$-group 
$\Delta_\rho^\phi$ with the $[\psi\phi^{-2},\widetilde{\nu}]$-isotypic part for
the action of  a prime to $p$ order group. This   geometric description 
of the Hecke action of $\Delta_\rho^\phi$ will play an important role 
for our construction.

\section{Modularity of the minimally ramified deformations.} \label{min}

Let $\rho: \G_F \rightarrow  \GL_2(\kappa)$ be a
continuous representation satisfying $\mathrm{{\bf (LI}_{Ind \rho}{\bf )}}$
and $\mathrm{{\bf (Mod}_{\rho}{\bf )}}$. 

The main aim of this section is to prove  :

\begin{theo}\label{minimal}
 Suppose     
$P_{\rho}=\varnothing$. Then all minimally ramified deformations
 of  $\rho$ are modular.
\end{theo}

In the notations of \S\ref{strategy}, the above theorem amounts to prove that
 $\pi: \cR \twoheadrightarrow \cT$ is  an isomorphism
(since $\Sigma=\varnothing$ in  the entire  section,  we shall omit
the subscripts). Our proof  uses a stronger version, due to Fujiwara \cite[\S2]{Fu},
of a method invented by Wiles \cite{wiles2} and
Taylor-Wiles \cite{TaWi} and known as a  Taylor-Wiles system 
(a similar formalism has been found independently by Diamond
\cite{diamond}). 

The construction of a Taylor-Wiles system will occupy the entire
section. It includes namely a geometric realization of  $\cT$ as a
Hecke algebra acting on the local component $\cM$ at
$\rho$ of the middle degree cohomology of a Hilbert
modular variety. The torsion freeness of $\cM$   is a crucial
ingredient (\cf Theorem \ref{theo-AB}(i)).
Lemmas \ref{aux}, \ref{isom1}, \ref{isom4} and \ref{isom5}   
are proved using standard fact about automorphic representations
and local Langlands correspondence for $\GL(2)$, whereas 
propositions \ref{isom3}, \ref{isom2} and \ref{isom6} use
finer geometric arguments.

\medskip
Note that Fujiwara's formalism isn't essential for us since 
we know that $\cM$ is free over $\cT_{P_{\rho}}$ and
$\cT_{P_{\rho}}$ is Gorenstein. This fact is an 
important ingredient in the proof of Theorem A, and is shown in
Proposition \ref{isom3} {\it without} assuming
$P_{\rho}=\varnothing$. Actually, we will only assume 
$P_{\rho}=\varnothing$ in \S\ref{selmer-groups}.

\medskip
\subsection{The formalism of Taylor-Wiles systems, following Fujiwara.}

\begin{defin} Let $\cQ$ be a family of finite sets of primes $\gq$ of  $F$
such that  $\Nm(\gq)\equiv 1\pmod{p}$.
A Taylor-Wiles system for  $\cQ$   is a family
$\{\cR,\cM,(\cR_Q,\cM^Q)_{Q\in\cQ}\}$ such that

(TW1)  $\cR_Q$ is a local complete $\cO[\Delta_Q]$-algebra, where
$\Delta_Q=\prod_{\gq  \in Q}\Delta_{\gq}$  and $\Delta_{\gq}$
is the $p$-Sylow of $(\go/\gq)^\times$.

(TW2) $\cR$ is a local complete $\cO$-algebra and there is an
isomorphism of local complete $\cO$-algebras 
$\cR_Q\otimes_{\cO[\Delta_Q]}\cO \cong \cR$.

(TW3) $\cM$ is a non-zero $\cR$-module, and $\cM^Q$ is an
$\cR_Q$-module, free of finite rank over $\cO[\Delta_Q]$ and such that 
$\cM^Q\otimes_{\cO[\Delta_Q]}\cO$ is isomorphic to $\cM$ as  $\cR$-module.

 We denote by $\cT$ the image of $\cR\rightarrow
\End_{\cO}(\cM)$.
\end{defin}

When $\cQ=\{Q_m|m\in \n\}$, we will write $\cR_m,\cM_m,...$ instead of
  $\cR_{Q_m},\cM^{Q_m},...$  .

\begin{theo} \cite[\S 2]{Fu}\label{stw}
Let $\{\cR,(\cR_m,\cM_m)_{m\in\n}\}$ be a Taylor-Wiles
  system. Assume that for all $m$:

\begin{enumerate}
\item  for all $\gq \in Q_m$, $\Nm(\gq)\equiv 1\pmod{p^m}$,
\item  $\cR_m$ can be  generated by $\#Q_m=r$ elements as a local complete
$\cO$-algebra.
\end{enumerate}

Then, the natural surjection $\cR\twoheadrightarrow \cT$ is an
isomorphism. Moreover, these algebras are flat  and complete
intersection of relative dimension zero over $\cO$ and  $\cM$ is
free over  $\cT$.
\end{theo}

\subsection{The rings $\cR_Q$.}\label{R_Q}

Let $Q$ be a finite set of auxiliary primes $\gq$ of $F$ satisfying  : 
\begin{enumerate}
\item  $\Nm(\gq)\equiv 1\pmod{p}$, and
\item  $\phi$ and $\rho$ are unramified at $\gq$, and
$\rho(\Frob_{\gq})$ has two distinct eigenvalues
$\alpha_{\gq}$ and $\beta_{\gq}$ in $\kappa$.
\end{enumerate}

\medskip
For such a $Q$ we can associate by \S\ref{strategy} an universal deformation ring
 $\cR_Q$,  endowed with a canonical surjection
$\cR_Q\twoheadrightarrow \cR_{\varnothing}=:\cR$. By a result of
 Faltings (\cf \cite[Appendix]{TaWi})
 $\cR_Q$ is a $\cO[\Delta_Q]$-algebra and
 $\cR_Q\otimes_{\cO[\Delta_Q]}\cO \cong \cR$. Thus (TW1) and  (TW2) hold.

More generally, for any set of primes $P$ disjoint from $Q$, 
$\cR_{Q\cup P}$ is a $\cO[\Delta_Q]$-algebra and
 \begin{equation}
\cR_{Q\cup P}\otimes_{\cO[\Delta_Q]}\cO \cong \cR_{P}.
\end{equation} 

In particular $\cR_{Q\cup P_{\rho}}$ is a $\cO[\Delta_Q]$-algebra.

\subsection{The module $\cM$.}

Denote by  $Y_{\rho}$ the  Hilbert modular varieties of level
$K_\rho$ defined in \S\ref{level}.

Let  $S$ be a finite set  of primes containing
$\Sigma_{\rho}\cup\{v\!\mid\! p\}\cup\{\gu\}\cup\{v\!\mid \phi\widetilde{\nu}_v
\text{ ramified}\}=\Sigma_{K_\rho}\cup\{v\!\mid\! p\}$.
Denote $\gm_{\rho}$ the maximal ideal of 
$\T^{S}=\cO[T_v, S_v \enspace | \enspace v\notin S]$
 corresponding to ${\rho}$. 

Recall that $[\psi ]$ denotes the  $\psi $-isotypic part for  the
 action of the Hecke operators $S_v$, $v\notin S$, and 
 $[\widetilde{\nu}\phi]$ the intersection of
 $[\phi\widetilde{\nu}_v]$-isotypic  parts for  the action of the Hecke 
operators $U_\delta=[K_\rho\left(\begin{smallmatrix}\delta & 0 \\ 0 & 1
\end{smallmatrix}\right)K_\rho]$ for $\delta\in \go_v^\times$.

\medskip

Fix an eigenvalue $\alpha_{\gu}$  of $\rho(\Frob_{\gu})$ and 
consider the $\cO$-module : 
\begin{equation}
\cM:=\rH^d(Y_{\rho},\bV_{\cO})[\psi ,\widetilde{\nu}\phi]_{
(\gm_\rho,U_{\gu}-\alpha_{\gu})}.
\end{equation}

Let $\T'$  be the image of  $\T^{S}$ in the ring of
$\cO$-linear endomorphisms of $\cM$.

By Theorem \ref{theo-AB}(i) the $\gm_{\rho}$-localization of the 
$\T^{S}$-module $\rH^d(Y_{\rho},\bV_{\cO})$ is free over $\cO$. 
Hence $\cM$ is free over $\cO$ as a direct factor of free $\cO$-module.

Moreover, $\cM$ is non-zero by  $\mathrm{{\bf (Mod}_{\rho}{\bf )}}$ 
and remark \ref{lowering}.   
For any newform   $f$ occurring in $\cM$,
consider the maximal ideal  
$$\gm_{f}=(\varpi, T'_{v}-\iota_p(c(f,v)),
S'_v-\iota_p(\psi(v))\N_{F/\Q}(v)^{1-k_0},U'_{v'}-\iota_p(c(f,v')) ; 
v\notin \Sigma_\rho, v'\in \Sigma_\rho)$$
of  $\T^{\full}=\cO[T'_{v},S'_{v} ;v\notin \Sigma_\rho][U'_{v'} ;v'\in
\Sigma_\rho]$. Note that  $\gm_{f} \cap \T^{S}= \gm_{\rho}$.

Let $\T$ (resp.  $\gm$) be the image of  
$\T^{\full}$ (resp.  $\gm_{f}$) in the ring of $\cO$-linear endomorphisms of 
$\rH^d(Y_{\rho},\bV_{\cO})[\psi ,\widetilde{\nu}\phi]$.

\begin{lemma}\label{isom1}
\begin{enumerate}
\item  There is an unique isomorphism 
of $\T^{S}$-algebras
$\cT_{P_{\rho}}\overset{\sim}{\rightarrow}\T'$,  
\item $\cM\otimes \C$ is free of rank $2^d$ over
  $\cT_{P_{\rho}}\otimes \C$, and 
\item the natural injective algebra homomorphism 
$\T' \hookrightarrow \T_{\gm}$ is an   isomorphism.
\end{enumerate}
\end{lemma}

\begin{proof} (i)
By lemma \ref{aux}, we have $\cT_{P_{\rho}\cup \{\gu\}}\cong\cT_{P_{\rho}}$.
Since  $\cO$-algebras $\cT_{P_{\rho}\cup \{\gu\}}$ 
and $\T'$ are torsion free (the first one by
definition,  the second one because $\cM$ is free over $\cO$), 
it is enough to show that  there is an unique isomorphism
 of  $\T^{S}\otimes\C$-algebras between 
$\cT_{P_{\rho}\cup \{\gu\}}\otimes\C$ and $\T'\otimes\C$
(tensors being over $\cO$ for some fixed embedding $\cO\hookrightarrow
\C$).

Consider a (cuspidal) automorphic representation $\pi$ generated by a 
holomorphic newform $f$ of weight $k$, central character
$\psi |\cdot|^{2-k_0}$ 
and prime to $p$ conductor. By definition $\pi$ contributes to  
$\cT_{P_{\rho}\cup \{\gu\}}\otimes\C$ if, and only if,
for all primes $v \nmid  py$, $v\notin P_{\rho}$, 
$\phi^{-1}\otimes\rho_{f,p}|_{\G_{F_v}}$ is a minimally ramified
deformation of $\rho_v$.

For $v\notin P_{\rho}$, $v \neq y$, remark \ref{mini-rk} shows that
$\phi^{-1}\otimes\rho_{f,p}|_{\G_{F_v}}$ is a minimally ramified 
deformation of $\rho_v$ if, and only if,
$(\phi\widetilde{\nu}_v)^{-1}\otimes\rho_{f,p}|_{\G_{F_v}}$  
has conductor $c_v$. By  Carayol's Theorem \cite{carayol} on the
compatibility between 
 the local and the global Langlands correspondences this is equivalent to 
$(\pi_v\otimes(\phi\widetilde{\nu}_v)^{-1})^{K_1(v^{c_v})}\cong
\pi_v^{K_v'}[\phi\widetilde{\nu}_v]\neq 0$.

If $v\in P_{\rho}$ then 
$\dim(\rho_{f,p}\otimes (\phi\widetilde{\nu}_v)^{-1})^{I_v}=\dim
(\rho\otimes \nu_v^{-1})^{I_v}=0$,  
hence $(\phi\widetilde{\nu}_v)^{-1}\otimes\rho_{f,p}|_{\G_{F_v}} $ has conductor $c_v$
and so $(\pi_v\otimes(\phi\widetilde{\nu}_v)^{-1})^{K_1(v^{c_v})}\neq 0$.

Finally, the argument of lemma \ref{aux} shows that $\pi_{\gu}$ is
unramified, hence $\pi_{\gu}^{K_0(\gu)}$ is two dimensional and 
contains an unique eigenline for $U_{\gu}$
with eigenvalue $\widetilde{\alpha}_{\gu}$ congruent to $\alpha_{\gu}$
modulo $\varpi$.  

Therefore,  $\pi$ contributes to $\cM \otimes\C$. 
By the Matsushima-Shimura-Harder  isomorphism, this is equivalent to 
$\pi$    contributing  to $\T'\otimes\C$. 
Conversely, if $\pi$    contributes to $\T'\otimes\C$, 
the same arguments show that  $\pi$    contributes to
$\cT_{P_{\rho}\cup \{\gu\}}\otimes\C$. 

(ii) 
Let $\pi$ be an automorphic representation contributing to  $\T'\otimes \C$. 
As a byproduct of the computations in (i) we have
 $\dim\pi_{\gu}^{K_0(\gu)}[U_{\gu}-\widetilde{\alpha}_{\gu}]=1$ 
and $\dim\pi_v^{K_v'}[\phi\widetilde{\nu}_v]=1$ for all $v \neq y$.
By the Matsushima-Shimura-Harder  isomorphism, the 
 $\pi_{\rf}$-part of $\cM \otimes\C$ is $2^d$-dimensional. 

(iii) 
We have to show that for all $v\in S$  the image 
$T'_v$ (or $U'_v$) in $\End_{\cO}(\cM)$ belong to $\T'$. 
The argument uses  local Langlands correspondence and the 
fact that $\cM$ is torsion free. 
As observed in \S\ref{strategy} there  exists a $P_{\rho}$-deformation
$\widetilde{\rho}$ of $\rho$ with coefficients in $\cT_{P_{\rho}}$ and by (i) 
there is an unique isomorphism of 
$\T^{S}$ algebras $\cT_{P_{\rho}}\cong \T'$.
It remains to prove that the resulting homomorphism 
$\cT_{P_{\rho}} \rightarrow \T_{\gm}$ is surjective.

If $v\notin \Sigma_\rho$, then the eigenvalue of  $T'_v$ on 
$\pi_v^{K_v'}[\phi\widetilde{\nu}_v]$ equals the eigenvalue of  $T_v$ on 
$(\pi_v\otimes(\phi\widetilde{\nu}_v)^{-1})^{K_1(v^{c_v})}$. Recall that 
$\nu_v(\varpi_v)=1$. 
Hence  the action of  $T'_v$ on $\cM$ is given by
$\tr(\widetilde{\rho}\otimes(\phi\widetilde{\nu}_v)^{-1})(\Frob_v)
\in \cT_{P_{\rho}} $.

If $v\in S_\rho$, then the eigenvalue of  $U'_v$ on 
$\pi_v^{K_v'}[\phi\widetilde{\nu}_v]$ equals the eigenvalue of  $U_v$ on 
$(\pi_v\otimes(\phi\widetilde{\nu}_v)^{-1})^{K_1(v^{c_v})}$.
Hence  the action of  $U'_v$ on $\cM$ is given by
the eigenvalue of 
$(\widetilde{\rho}\otimes(\phi\widetilde{\nu}_v)^{-1})(\Frob_v)$  on 
the line  $(\widetilde{\rho}\otimes(\phi\widetilde{\nu}_v)^{-1})^{I_v}$
hence belongs to $\cT_{P_{\rho}}$.

If $v\in  \Sigma_\rho\backslash S_\rho$, then $U'_v=0$. This completes
the proof. 
\end{proof}

\begin{prop} \label{isom3} 
$\cM$ is free of rank $2^d$ over $\cT_{P_{\rho}}$
 and $\cT_{P_{\rho}}$ is Gorenstein.
\end{prop}

\begin{proof} Put  
$W=\rH^d(Y_{\rho}, \bV_{\kappa})[\psi ,\widetilde{\nu}\phi]_{(\gm_{\rho},U_{\gu}-\alpha_{\gu})}$. 
By lemma \ref{isom1} and \cite[lemma 6.8]{dim-hmv}, it is enough to show that 
$W[\gm]=\cM\otimes_{\T_{\gm}}\kappa$ is a $\kappa$-vector space of
dimension at most $2^d$.  

As in the proof of Proposition \ref{ihara}, the condition
$\mathrm{\bf (LI_{\Ind \rho})}$ implies that every
$\G_{\widetilde{F}}$-irreducible subquotient 
of $W[\gm]\subset W[\gm_{\rho}]$ is  isomorphic to $\otimes \Ind_F^{\Q}\rho$.
Therefore it is enough to check that the last graded piece 
of the  Fontaine-Laffaille module attached to  $W[\gm]$ 
has dimension $\leq 1$. Again as in the proofs of Theorem \ref{theo-AB}
and Proposition \ref{ihara}, this amounts to showing that :
\begin{equation}
\dim
\rH^0(\cY_\rho{}_{/\kappa},\underline{\omega}^k\otimes\underline{\nu}^{1-k_0/2})[\psi,\nu,\gm]\leq
1.\end{equation} 

By the $q$-expansion Principle, a Hilbert modular form 
in $\rH^0(\cY_\rho{}_{/\kappa},\underline{\omega}^k\otimes\underline{\nu}^{1-k_0/2})$
 is uniquely determined by the coefficients
of its $q$-expansion. The coefficients are indexed by 
$(F\otimes \widehat{\Z})/ \prod_v \ker(\nu_v)$, hence 
a form in $\rH^0(\cY_\rho{}_{/\kappa},\underline{\omega}^k\otimes\underline{\nu}^{1-k_0/2})[\nu]$
is uniquely determined by the subset of its coefficients
indexed by  ideal of $F$, and is it a standard fact that 
coefficients at non-integral ideals vanish. 

Finally, the coefficients of a form in 
$\rH^0(\cY_\rho{}_{/\kappa},\underline{\omega}^k\otimes\underline{\nu}^{1-k_0/2})[\nu][\psi,\gm]$
are uniquely determined, since they are related to the 
eigenvalues of $T'_v$, $S'_v$ and $U'_v$, and those are fixed in
the $[\psi,\gm]$-part. 
\end{proof}

\subsection{The modules $\cM^Q$.}

Denote by $Y_{0,Q}$  (resp.  $Y^Q$)  the  Hilbert modular
varieties of level $K_{0,Q}$ (resp.   $K^Q$) introduced in \S\ref{level}.
The natural homomorphism  $Y^Q\rightarrow Y_{0,Q}$ induced
by the inclusion $K^Q\subset K_{0,Q}$, is étale with  group  $\Delta_Q$.

\medskip

Assume    that  $S$ contains
$\Sigma_{\rho}\cup\{v\!\mid\! p\}\cup Q\cup\{\gu\}\cup
\{v\!\mid \phi\widetilde{\nu}_v
\text{ ramified}\}=\Sigma_{K^Q}\cup\{v\!\mid\! p\}$.

\medskip
Let $\T'_{0,Q}$   be the image of the Hecke algebra $\T^{S}$
in the ring of $\cO$-linear endomorphisms of  :
\begin{equation}
\cM_{0,Q}:=\rH^d(Y_{0,Q},\bV_{\cO})[\psi ,\widetilde{\nu}\phi]_{
(\gm_\rho,U_{\gu}-\alpha_{\gu},U_{\gq}-\alpha_{\gq} ;
\gq\in Q)}.\end{equation}

Let $\T'_Q$ be the image of the Hecke algebra $\T^{S}[\Delta_Q]$
in the ring of $\cO$-linear endomorphisms of
 \begin{equation}
\cM^Q:=\rH^d(Y^Q,\bV_{\cO})[\psi ,\widetilde{\nu}\phi]_{(\gm_\rho,
   U_{\gu}-\alpha_{\gu},U_{\gq}-\alpha_{\gq} ; \gq\in Q)}.\end{equation}

Note that the group $\Delta_{Q}$ acts on 
$\rH^d(Y^Q,\bV_{\cO})$ via the Hecke operators
$U_\delta$, $\delta\in \gq$, $\gq\in Q$ defined in \S\ref{level}

Note that whereas $U_{\gq}\in
\End_{\cO}(\rH^d(Y^Q,\bV_{\cO})[\psi ,\widetilde{\nu}\phi]_{\gm_\rho})$  
depends on the choice of an uniformizer, the 
ideal $(\varpi,U_{\gq}-\alpha_{\gq})$ doesn't, so $\cM^Q$ doesn't.

Again by  Theorem \ref{theo-AB}(i) the modules $\cM_{0,Q}$ and $\cM^Q$
are  free over $\cO$, hence  $\T'_{0,Q}$ and  $\T'_Q$ are torsion free.

By lemma \ref{isom1}, for all $\gq\in Q$, the Hecke operators 
$T_{\gq}$ and $S_{\gq}$
belong to $\cT_{P_\rho}\overset{\sim}{\rightarrow}\T'$, hence act on $\cM$.
By  \S\ref{R_Q} and Hensel's lemma the polynomial
$X^2-T_{\gq}X+S_{\gq}\Nm(\gq)\in \cT_{P_\rho}[X]$ has
an unique root $\widetilde{\alpha}_{\gq}\in \cT_{P_\rho}$
(resp. $\widetilde{\beta}_{\gq}\in \cT_{P_\rho}$) above $\alpha_{\gq}$ (resp. $\beta_{\gq}$).

\begin{lemma} \label{isom4}
There  exists  an unique isomorphism 
of $\T^{S}$-algebras
$\T'_{0,Q}\overset{\sim}{\rightarrow} \T'$.
\end{lemma}

\begin{proof}
As in lemma \ref{isom1}(i) it is enough to show that there is an  isomorphism
 of $\T^{S}$-algebras $\T'_{0,Q}\otimes\C \overset{\sim}{\rightarrow} 
\T'\otimes\C$. 

The local component at $\gq$ 
of an automorphic representation $\pi$  contributing to $\T'_{0,Q}\otimes\C$
(or $\cM_{0,\gq}\otimes \C$)  admits invariants by $K_0(\gq)$ and
cannot be special  (since  $\alpha_{\gq}\neq
\beta_{\gq}\Nm(\gq)^{\pm 1}$ by our assumptions in \S\ref{R_Q}); hence
it is necessarily  
an unramified principal series and so contributes to  $\cM\otimes \C$
and $\T'\otimes\C$. 
Moreover, $\pi$ contributes with the same multiplicity  
both in $\cM_{0,\gq}\otimes \C$ and  $\cM\otimes \C$. 
The proof of this fact is very similar to the proof of lemma 
\ref{isom1}(ii), once we notice that for every such $\pi$, 
 $\pi_{\gq}^{K_0(\gq)}$ is two dimensional and 
contains an unique eigenline for $U_{\gq}$
with eigenvalue congruent to $\alpha_{\gq}$
modulo $\varpi$. 
\end{proof}

\begin{lemma}\label{isom5}
There is an unique isomorphism 
of  $\T^{S}[\Delta_Q]$-algebras
$\cT_{P_{\rho} \cup Q}\overset{\sim}{\rightarrow}\T'_Q$.
\end{lemma}

\begin{proof} Both $\cT_{P_{\rho} \cup Q}$ and $\T'_Q$
are defined as images of $\T^{S}[\Delta_Q]$ hence the uniqueness. 
 For the existence, as in lemma \ref{isom1}(i), 
 it is enough to show that there is an  isomorphism
 of  $\T^{S}[\Delta_Q]$-algebras between 
$\cT_{P_{\rho} \cup Q}\otimes\C$ and $\T'_Q\otimes\C$.

Consider a (cuspidal) automorphic representation $\pi$ generated by a 
holomorphic newform $f$ of weight $k$, central character
$\psi |\cdot|^{2-k_0}$ and prime to $p$ conductor.

If $\pi$ contributes $\T'_Q\otimes\C$ then it 
necessarily contributes to $\cT_{P_{\rho} \cup Q}\otimes\C$, 
since by the proof of lemma \ref{isom1}(i)  $\rho_{f,p}$
satisfies all the deformation conditions at primes outside $Q$, 
and there is no  deformation conditions at primes in $Q$.

Conversely, suppose that  $\pi$ contributes to 
$\cT_{P_{\rho} \cup Q}\otimes\C$. By  \cite[Appendix]{TaWi},
$\rho_{f,p}|_{\G_{F_{\gq}}}$ is decomposable 
and   $\rho_{f,p}|_{I_{\gq}}\cong\chi\oplus\chi^{-1}$ where $\chi$ factors
through the natural surjective homomorphism $I_{\gq}\rightarrow\go_{\gq}^\times
\rightarrow (\go/\gq)^\times\rightarrow \Delta_{\gq}$.
By the local Langlands correspondence $\pi_{\gq}$ is
a principal series induced from two characters whose restriction 
to $\go_{\gq}^\times$ are $\chi$ and $\chi^{-1}$. It follows that 
\begin{equation}
\pi_{\gq}^{K_{\gq}}=\begin{cases}\pi_{\gq}^{K_0(\gq)}\text{ , if }
\chi \text{ is trivial, and }\\
(\pi_{\gq}\otimes \chi)^{K_1(\gq)}\oplus 
(\pi_{\gq}\otimes \chi^{-1})^{K_1(\gq)}\text{ , if }
\chi \text{ is non-trivial.}\end{cases}\end{equation}
In both cases $\pi_{\gq}^{K_{\gq}}$ is two dimensional and splits
under the action of $U_{\gq}$ as a direct sum of two lines, 
one with eigenvalue $\widetilde{\alpha}_{\gq}$ congruent to
$\alpha_{\gq}$ modulo $\varpi$  and 
one with eigenvalue  $\widetilde{\beta}_{\gq}$ congruent to
$\beta_{\gq}$ modulo $\varpi$.  Hence 
$\pi_{\gq}^{K_{\gq}}[U_{\gq}-\alpha_{\gq}]\neq 0$. 
Note that whereas $U_{\gq}$ and the eigenvalue depend
on the choice of an uniformizer, the decomposition doesn't. 

Also, note that by local Langlands correspondence, 
the $\Delta_{\gq}$-action on $\T'_Q\otimes\C$
coming from the Hecke action of $K_0(\gq)$ on 
$\pi_{\gq}^{K_{\gq}}$,  corresponds to the $\Delta_{\gq}$-action 
on $\cT_{P_{\rho} \cup Q}\otimes\C$ coming from the 
 $I_{\gq}$-action on $\rho_{f,p}$. 

The above discussion at primes in $Q$ together with the arguments
of lemma \ref{isom1}(i) at the primes outside $Q$ imply that 
 $\pi$ contributes to $\cM^Q\otimes\C$, hence to $\T'_Q\otimes\C$.

\end{proof}

\subsection{The condition (TW3).}

\begin{prop} \label{isom2}
There is a $\T^{S}$-linear isomorphism 
$\cM\overset{\sim}{\rightarrow}\cM_{0,Q}$ such that the 
$U_{\gq}$-action on $\cM_{0,Q}$ correspond to the  
$\widetilde{\alpha}_{\gq}$-action on $\cM$. 
\end{prop}

\begin{proof} We may assume that $Q=\{\gq\}$ and prove the lemma
with $K_\rho$ replaced by 
$K_{0,Q\bs \{\gq\}}$ in the definitions of $Y_{\rho}$,
 $\T'$ and $\cM$. Consider the $\T^{S}$-linear homomorphism :
$$ \cM \rightarrow \cM^2\textrm{  , }
x \mapsto (x, -x|\widetilde{\beta}_{\gq}).$$

Let  $U_{\gq}$ be the  $\T^{S}$-linear endomorphism of $\cM^2$ given 
by the matrix $\left(\begin{smallmatrix} T_{\gq}  & - \Nm(\gq) S_{\gq}
\\ 1 & 0\end{smallmatrix}\right)$ acting on the right. 
 Since its eigenvalues  $\widetilde{\alpha}_{\gq}$ and 
 $\widetilde{\beta}_{\gq}$  are distinct modulo $\varpi$, it
 induces an isomorphism :
$$\cM \overset{\sim}{\rightarrow}  (\cM^2)_{(U_{\gq}-\alpha_{\gq})}.$$

Consider the natural degeneracy maps  
$\pr_1,\pr_2 : Y_{0,\gq}\rightarrow Y_{\rho}$ used in the 
definition of the Hecke correspondence $T_{\gq}$ in \S\ref{hmv}.
The  $\T^{S}$-linear homomorphism
$\pr_1^*+\pr_2^*: \rH^d(Y_{\rho},\bV_{\cO})^2 \rightarrow
\rH^d(Y_{0,\gq},\bV_{\cO})$ yields 
(after taking $[\psi ,\widetilde{\nu}\phi]$-parts and  localizing  at $\gm_{\rho}$) :
$$\xi : \rH^d(Y_{\rho},\bV_{\cO})[\psi ,\widetilde{\nu}\phi]_{\gm_{\rho}}^2\rightarrow 
\rH^d(Y_{0,\gq},\bV_{\cO})[\psi ,\widetilde{\nu}\phi]_{\gm_{\rho}}.$$

From the definition of $U_{\gq}$ acting on  $\cM^2$ we see that  $\xi$
is $U_{\gq}$-linear. It is also  $U_{\gu}$-linear, hence after 
localization at $(\varpi, U_{\gq}-\alpha_{\gq},U_{\gu}-\alpha_{\gu})$
  induces :
$$\xi' : (\cM^2)_{(U_{\gq}-\alpha_{\gq})}
\rightarrow \cM_{0,\gq}.$$
It is enough to  show then that  $\xi'$  is an isomorphism.
\smallskip

By lemma \ref{isom4} and its proof, we see that  $\xi'\otimes \id_{\C}$ is an 
isomorphism. It remains to prove that $\xi$ (hence $\xi'$) 
is injective with flat cokernel.

Let $\widehat{\xi}$ be the dual of $\xi$
with respect to the modified Poincar{\'e} 
pairing defined in \S\ref{poincare}. The matrix of
$\widehat{\xi}\circ\xi 
: (\cM\otimes\kappa)^2\rightarrow
(\cM\otimes\kappa)^2$  is given by
$\left(\begin{smallmatrix}1+\Nm(\gq) & S_{\gq}^{-1}T_{\gq}
\\T_{\gq} &1+\Nm(\gq)  \end{smallmatrix}\right)$. It is invertible by our
assumptions on $\gq$. Therefore $\xi$ is injective with flat cokernel. 
\end{proof}

By \S\ref{R_Q}, $\cR_{P_{\rho} \cup Q}$ is a $\cO[\Delta_Q]$-algebra.
Hence the  surjective  homomorphism of  local $\cO$-algebras
$\pi_{\Sigma}:\cR_{P_{\rho} \cup Q}\rightarrow \cT_{P_{\rho} \cup Q}$ defined
in  \S\ref{strategy} endows  $\cT_{P_{\rho} \cup Q}$ with
$\cO[\Delta_Q]$-algebra structure. 

\begin{prop}  \label{isom6} $\cM^Q$ is a free $\cO[\Delta_Q]$-module
and $\cM^Q\otimes_{\cO[\Delta_Q]}\cO\cong \cM_{0,Q}$ as   $\T^{S}$-modules.
\end{prop}

\begin{proof}
By Theorem \ref{freeness}(i) $\rH^d(Y^Q,\bV_{\cO})_{\gm_{\rho}}$ 
is free over $\cO[\Delta_Q]$ and the  $\T^{S}$-module of its
$\Delta_Q$-coinvariants is  isomorphic to
$\rH^{d}(Y_{0,Q},\bV_{\cO})_{\gm_{\rho}}$. If  the class group
$\Cl_\rho$ defined in \S\ref{twisted}  has order prime to $p$
(in particular $\phi$ is trivial) then the claim  follows simply by
taking the $[\psi,\widetilde{\nu}]$-part. In fact  the
$[\psi,\widetilde{\nu}]$-part, for the action of a prime to $p$ order group, 
of a  free $\cO[\Delta_Q]$-module is a free $\cO[\Delta_Q]$-direct factor.

In the general case, denote by  $\Delta_\rho^\phi$ the $p$-Sylow subgroup 
of $\Cl_\rho\times \prod_{v\in \Sigma_\rho} (\go_v^\times/\ker(\phi_v))$.
As in \S\ref{twisted} the $p$-group 
$\prod_{v\in \Sigma_\rho} (\go_v^\times/\ker(\phi_v))$ injects 
in $\A^\times K_0(Q\gn_{\varnothing})/\A^\times K_{0,Q}$
and {\it a fortiori} in $\A^\times K_0(Q\gn_{\varnothing})/\A^\times K^Q$.
Also the  morphisms $Y_{Q}\rightarrow Y_{Q}^{\ad}$ and 
 $Y_{0,Q}\rightarrow Y_{0,Q}^{\ad}$ are étale with group
 $\Cl_\rho$. Hence  the étale morphism
$Y_{Q}\rightarrow Y_{Q\gn_{\varnothing}}^{\ad}$ (resp.  
$Y_{0,Q}\rightarrow Y_{Q\gn_{\varnothing}}^{\ad}$)  
factors through an étale morphism $Y_{Q}\rightarrow Y_{Q}^{\Delta}$
(resp. $Y_{0,Q}\rightarrow Y_{0,Q}^{\Delta}$) with group  $\Delta_\rho^\phi$.
 Then  Theorem \ref{freeness}(i) applies to 
each of the five étale morphisms in the following diagram : 
\begin{equation}
\xymatrix@C=20pt@R=8pt{
& Y_{Q}\ar[dl]_{\Delta_Q}\ar[dr]^{\Delta_\rho^\phi}
\ar[dd]& \\
Y_{0,Q}\ar[dr]_{\Delta_\rho^\phi}& & Y_{Q}^{\Delta} \ar[dl]^{\Delta_Q} \\&
Y_{0,Q}^{\Delta} & }\end{equation} 
In particular,  $\rH^d(Y^Q,\bV_{\cO})_{\gm_{\rho}}$ is 
free over $\cO[\Delta_\rho^\phi\times\Delta_Q]$, hence
$\rH^d(Y^Q,\bV_{\cO})_{\gm_{\rho}}[\phi]$ is 
free over $\cO[\Delta_Q]$ and 
$$\rH^d(Y^Q,\bV_{\cO})_{\gm_{\rho}}[\phi]
\otimes_{\cO[\Delta_Q]}\cO \cong \left(\rH^d(Y^Q,\bV_{\cO})_{\gm_{\rho}}
\otimes_{\cO[\Delta_Q]}\cO\right)[\phi]
\cong \rH^d(Y_{0,Q},\bV_{\cO})_{\gm_{\rho}}[\phi].$$

Taking further the $[\psi\phi^{-2},\widetilde{\nu}]$-part, for the action of the 
prime to $p$ order group
$(\Cl_\rho/\Delta_\rho^\phi)\times \prod_{v}
(\go_v^\times/\ker(\nu_v))$, and using the argument 
invoked in the beginning of the proof, yields the desired result.
\end{proof}

\medskip
So far we have constructed a Taylor-Wiles system 
$\{\cR,\cM,(\cR_Q,\cM^Q)_{Q\in\cQ}\}$ for the family $\cQ$ of
sets $Q$ containing a finite number of primes $\gq$ as in
\S\ref{R_Q}. The aim of the next  paragraph is to find a subfamily
$\{Q_m|m\in \n\}$ satisfying the conditions  (i) and (ii) of Theorem \ref{stw}.

\subsection{Selmer groups.} \label{selmer-groups}

We assume in this paragraph that $P_{\rho}=\varnothing$. 
Let  $\rho_{f,p}$ be a modular deformation of $\rho$ as in 
$\mathrm{{\bf (Mod}_{\rho}{\bf )}}$. 
For $r\geq 1$ we put $\rho_r := \rho_{f,p} \mod{\varpi^r}$,  so
that $\rho_1=\rho$.

We will use  Galois cohomology techniques in order  to control
the number of generators of $\cR_Q$. 

\begin{defin}\label{finite}
For $v\!\mid\! p$ the subgroup
$\rH^1_{\rf}(F_v,\Ad^0\rho_r)\subset
\rH^1(F_v,\Ad^0\rho_r)$ consists of classes
corresponding to crystalline extensions of $\rho_r$ by itself.

For $v\nmid p$ the subgroup of unramified classes
$\rH^1_{\rf}(F_v,\Ad^0\rho_r)\subset
\rH^1(F_v,\Ad^0\rho_r)$  is defined as
$\rH^1(\G_{F_v}/I_v,(\Ad^0\rho_r)^{I_v})$.
\end{defin}

\begin{defin}\label{selmer}
The Selmer groups associated to a finite set of primes   $\Sigma$ are defined as
$$\rH^1_{\Sigma}(F,\Ad^0\rho_r)=
\ker\left(\rH^1(F,\Ad^0\rho_r) \rightarrow
  \bigoplus_{v\notin\Sigma}  \rH^1(F_v,\Ad^0\rho_r)/
\rH^1_{\rf}(F_v,\Ad^0\rho_r) \right),$$

$$\text{ and }\rH^1_{\Sigma}(F,\Ad^0\rho_{f,p}\otimes \Q_p/\Z_p)=\limind \rH^1_{\Sigma}(F,\Ad^0\rho_r).$$
\end{defin}

The dual of $\Ad^0\rho$ is canonically isomorphic to its
Tate twist  $\Ad^0\rho(1)$. The corresponding dual Selmer group
$\rH^1_{\Sigma^*}(F,\Ad^0\rho(1))$ is defined as the kernel
of the map
$$\rH^1(F,\Ad^0\rho(1)) \rightarrow
  \bigoplus_{v\in\Sigma} \rH^1(F_v,\Ad^0\rho(1))
\bigoplus_{v\notin\Sigma}  \rH^1(F_v,\Ad^0\rho(1))/
\rH^1_{\rf}(F_v,\Ad^0\rho(1)).$$

The Poitou-Tate exact sequence yields the following formula :
\begin{equation}
\frac{\#\rH^1_{\Sigma}(F,\Ad^0\rho)}{\#\rH^1_{\Sigma^*}(F,\Ad^0\rho(1))}=
\frac{\#\rH^0(F,\Ad^0\rho)}{\#\rH^0(F,\Ad^0\rho(1))}
\prod_{v\in \Sigma} \frac{\#\rH^1(F_v,\Ad^0\rho_v)}
{\#\rH^0(F_v,\Ad^0\rho_v)} \prod_{v\mid p\infty}
\frac{\#\rH^1_{\rf}(F_v,\Ad^0\rho_v)}
{\#\rH^0(F_v,\Ad^0\rho_v)}.\end{equation}
A proof for $F=\Q$ can be found in \cite[Proposition 1.6]{wiles2}, but
as mentioned in \cite[Theorem 2.19]{DDT}
the same argument works over an arbitrary number field.

 By  $\mathrm{\bf (LI_{\Ind \rho})}$ we have
$\rH^0(F,\Ad^0\rho)=\rH^0(F,\Ad^0\rho(1))=0$.
Since $\rho$ is totally odd, for all $v\mid \infty $ we
have $\dim \rH^0(F_v,\Ad^0\rho_v)=1$. Since
$\rho$ is  crystalline at all places $v$ dividing $p$ we have 
\begin{equation}
\dim\rH^1_{\rf}(F_v,\Ad^0\rho_v)
-\dim\rH^0(F_v,\Ad^0\rho_v)\leq [F_v:\Q_p]\end{equation} 
(\cf \cite[Theorem 3.20]{Fu} and also \cite[Cor.2.3]{dfg}). Finally,
for all $\gq\in Q$,  $\dim
\rH^0(F_{\gq},\Ad^0\rho_{\gq}(1))=1$.
Putting all together we obtain :

\begin{lemma}\label{estimate}
$\dim \rH^1_{Q}(F,\Ad^0\rho)\leq
\rH^1_{Q^*}(F,\Ad^0\rho(1))+\#Q.$
\end{lemma}

Finally, by the same arguments as in  \cite[\S3]{wiles2} we obtain : 

\begin{lemma}  Let  $m\geq 1$ be an integer. Then
for each non-zero  element $x\in
\rH^0_{\varnothing^*}(F,\Ad^0\rho(1))$ there exists a  prime  $\gq$
such that  :  

$\relbar$ $\Nm(\gq) \equiv 1\pmod{p^m}$,

$\relbar$ $\rho$ is unramified at $\gq$ and
$\rho(\Frob_{\gq})$ has two distinct eigenvalues in $\kappa$, and

$\relbar$ the image by the restriction map  of $x$ in
$\rH^1_{\rf}(F_{\gq},\Ad^0\rho(1))$ is  non-trivial.
\end{lemma}

Put $r:= \dim \rH^0_{\varnothing^*}(F,\Ad^0\rho(1))$.
For each $m\geq 1$, let $Q_m$ be the set of primes $\gq$ corresponding
by the above lemma to the elements of a basis of
$\rH^0_{\varnothing^*}(F,\Ad^0\rho(1))$.
Then $\rH^0_{Q_m^*}(F,\Ad^0\rho(1))=0$ and by lemma \ref{estimate}
we obtain $\dim \rH^0_{Q_m}(F,\Ad^0\rho)\leq\#Q_m$.
Therefore $\cR_m$ is generated by at most $\#Q_m=r$ elements.
This completes the proof of Theorem \ref{minimal}.

\section{Raising the level.}

\subsection{Numerical invariants.}

\begin{defin}\label{inv}
For a local complete noetherian $\cO$-algebra   $A$   endowed  with a
 surjective homomorphism  $\theta_A:A\rightarrow \cO$, we define  the
 following two      invariants :

$\relbar$ the congruence ideal
$\eta_A:=\theta_A(\Ann_A(\ker\theta_A))\subset \cO$, and

$\relbar$ the module of relative differentials 
$\Phi_A:=\Omega^1_{A/\cO}=\ker\theta_A/(\ker\theta_A)^2$.
\end{defin}

Here we state Wiles' numerical criterion :

\begin{theo} \label{numerical} \cite[Theorem 3.40]{DDT}
Let $\pi: \cR\twoheadrightarrow \cT$  be a surjective homomorphism
such that $\theta_{\cR}=\pi\circ\theta_{\cT}$. Assume that
 $\cT$ is finite and flat over $\cO$ and $\eta_{\cT}\neq (0)$.
Then the following three conditions are equivalent :

\begin{enumerate}
\item  $\#\Phi_{\cR}\leq\#(\cO/\eta_{\cT}) $,

\item  $\#\Phi_{\cR}=\#(\cO/\eta_{\cT}) $, and

\item  $ \cR$ and $\cT$ are complete
intersections over $\cO$ and $\pi$ is an isomorphism.
\end{enumerate}
\end{theo}

We consider couples $(\cT,\cM)$ consisting of a finite and  flat
$\cO$-algebra $\cT$ and a $\cT$-module $\cM$ which is a finitely
generated  free  $\cO$-module endowed with a perfect $\cT$-linear
pairing $\langle\cdot,\cdot\rangle:\cM \times \cM \rightarrow \cO $
and such that $\cM\otimes E$ is free over $\cT \otimes E$ of a given 
rank (in our application this rank will be $2^d$). The
pairing induces an isomorphism of $\cT$-modules
$\cM\overset{\sim}{\longrightarrow} \Hom(\cM,\cO)$.

\medskip
From \cite[Lemma 4.17]{DDT} and \cite[Theorem 2.4]{diamond} we deduce the
following :

\begin{prop} \label{raising}
Let  $(\cT,\cM)$ and  $(\cT',\cM')$ be two couples as above.
Assume that we have a surjective homomorphism $\cT'\twoheadrightarrow
\cT$ and a  $\cT'$-linear injective homomorphism $\xi:\cM\hookrightarrow
\cM'$ inducing via $\langle\cdot,\cdot\rangle$
a surjective homomorphism  $\widehat{\xi}:\cM' \twoheadrightarrow \cM$.

If  $\cM$ is free over $\cT$  and if
$\widehat{\xi}\circ\xi(\cM)=T\cdot\cM$ for some $T\in \cT$
then
$$\#(\cO/\eta_{\cT}) \#(\cO/\theta_{\cT}(T)) \leq
\#(\cO/\eta_{\cT'}).$$

Moreover,  equality holds if, and only if, $\cM'$ is free over $\cT'$.
\end{prop}

\subsection{Proof of theorem A.}\label{ideal}
 Let $\Sigma$ be a finite set of primes containing $P_{\rho}$. 
We start by  redefining  $\cT_{\Sigma}$  geometrically.

Let $Y_{\Sigma}$ be the Hilbert modular variety of level
$K_{\Sigma}$ defined in \S\ref{level}.

Let  $S$ be a finite set  of primes containing
$\Sigma_{\rho}\cup\Sigma\cup\{v\!\mid\! p\}\cup\{\gu\}\cup\{v\!\mid
\phi\widetilde{\nu}_v 
\text{ ramified}\}=\Sigma_{K_\Sigma}\cup\{v\!\mid\! p\}$.

\medskip
 Let   $\T'_{\Sigma}$ be the image of $\T^{S}$ 
in the ring of $\cO$-linear endomorphisms of  :
\begin{equation}
\cM_{\Sigma}:=\rH^d(Y_{\Sigma},\bV_{\cO})[\psi ,\widetilde{\nu}\phi]_{
(\gm_{\rho},U_{\gu}-\alpha_{\gu},U''_{\gq}; \gq\in \Sigma)}.
\end{equation}

By Theorem \ref{theo-AB}(i) $\cM_{\Sigma}$ is 
free of finite rank over $\cO$.

\medskip

For every    Hilbert modular newform $f$ occurring in $\cT_{\Sigma}$ 
 we denote by  $\theta_{f}^{\Sigma} : \cT_{\Sigma} \rightarrow \cO$
the projection on the $f$-component and by $\eta_{f}^{\Sigma}$ the  
corresponding congruence ideal.

\begin{lemma} \label{oldform}
\begin{enumerate}
\item There is an unique isomorphism of $\T^{S}$-algebras
$\T'_{\Sigma}\cong \cT_{\Sigma}$.

\item  $\cM_{\Sigma}\otimes \C$ is free of rank $2^d$ over 
$\cT_{\Sigma} \otimes \C$ and $U''_{\gq}$ acts as $0$ on it for all $\gq\in \Sigma$. 
\end{enumerate}
\end{lemma}

\begin{proof}
We follow closely the proofs of lemmas \ref{isom1} and \ref{isom5}. 
The main point here is to show that, if $f$ is a  Hilbert modular newform
occurring in $\cT_{\Sigma} \otimes \C$ and $\pi$ denotes the 
corresponding  automorphic representation, then 
then for all $\gq\in \Sigma$, 
$(\pi_{\gq})^{K''_{\gq}}[\phi_{\gq}\widetilde{\nu}_{\gq}]=
(\pi_{\gq}\otimes \phi_{\gq}^{-1}\widetilde{\nu}_{\gq}^{-1})^{K_1({\gq}^{c_{\gq}})
\cap K_0({\gq}^{c_{\gq}+d_{\gq}})}$ contains an unique eigenline for 
$U''_{\gq}$ with eigenvalue congruent to $0$ modulo $\varpi$ 
(and this eigenvalue is actually $0$). We distinguish three cases.

$\bullet$ If $(\widetilde{\nu}_{\gq}\phi_{\gq})^{-1}\otimes \rho_{f,p}$ is 
unramified at  $\gq$, then necessarily $d_{\gq}=2$, $c_{\gq}=0$ and 
$$\dim \left((\pi_{\gq}\otimes\phi_{\gq}^{-1}\widetilde{\nu}_{\gq}^{-1})^{K_0({\gq}^2)}\right)=3.$$
 The characteristic polynomial 
of  $U''_{\gq}=[K_0({\gq}^2)\left(\begin{smallmatrix} \varpi_{\gq} & 0 \\
0 & 1 \end{smallmatrix}\right) K_0({\gq}^2)]$ acting on it  is 
given by :
$$ X(X^2-c(f,\gq)X+\psi(\gq)\Nm(\gq)^{k_0-1})=\theta_f^\Sigma(
X(X^2-T'_{\gq}X+S'_{\gq}\Nm(\gq))),$$  
 and  $X=0$ is  simple  root modulo $\varpi$ of this polynomial.

$\bullet$ If  $\dim\left((\widetilde{\nu}_{\gq}\phi_{\gq})^{-1}
  \otimes\rho_{f,p}\right)^{I_{\gq}}=1$, 
then $d_{\gq}\geq 1$ and 
$$\dim \left((\pi_{\gq}\otimes\phi_{\gq}^{-1}\widetilde{\nu}_{\gq}^{-1})^{ 
K_1({\gq}^{c_{\gq}})\cap K_0({\gq}^{c_{\gq}+d_{\gq})}}\right)=2.$$
 The characteristic polynomial 
of  $U''_{\gq}=[K_0({\gq}^{c_{\gq}+d_{\gq}})
\left(\begin{smallmatrix} \varpi_{\gq} & 0 \\
0 & 1 \end{smallmatrix}\right)  K_0({\gq}^{c_{\gq}+d_{\gq}})]$  acting on it is 
given by : 
$$ X(X-c(f,\gq))=\theta_f^\Sigma(X(X-U'_{\gq})),$$  where 
$U'_{\gq}=[K_0({\gq}^{c_{\gq}+d_{\gq}-1})\left(\begin{smallmatrix} \varpi_{\gq} & 0 \\
0 & 1 \end{smallmatrix}\right) K_0({\gq}^{c_{\gq}+d_{\gq}-1})]$
 and  $X=0$ is  simple  root modulo $\varpi$ of this polynomial.

$\bullet$ Finally, if $\left((\widetilde{\nu}_{\gq}\phi_{\gq})^{-1}\otimes
\rho_{f,p}\right)^{I_{\gq}}=\{0\}$, then 
$$\dim \left((\pi_{\gq}\otimes\phi_{\gq}^{-1}
\widetilde{\nu}_{\gq}^{-1})^{K_1(\gq^{c_{\gq}}) 
\cap K_0(\gq^{c_{\gq}+d_{\gq}})}\right)=1,$$  and
$U''_{\gq}=0$ on it.  This completes the proof. 
\end{proof}

 By \S\ref{strategy}
we have a surjection $\pi_{\Sigma} : \cR_{\Sigma} \rightarrow
\cT_{\Sigma}$. Therefore we may endow $\cR_{\Sigma}$ with a
surjective homomorphism $\theta_{f}^{\Sigma}\circ\pi_{\Sigma} :
\cR_{\Sigma} \rightarrow  \cO$ and we denote $\Phi_{f}^{\Sigma}$
the corresponding numerical invariant.

\begin{prop}  \label{tangent}\cite[Proposition 1.2]{wiles2}
$\,\,\Hom_{\cO}(\Phi_{f}^{\Sigma},E/\cO)\cong
\rH^1_{\Sigma}(F,\Ad^0(\rho_{f,p})\otimes \Q_p/\Z_p).$
\end{prop}

Finally, one easily deduce from \S\ref{poincare} a  perfect
$\cT_{\Sigma}$-linear pairing :  
\begin{equation}
\langle\cdot,\cdot\rangle:\cM_{\Sigma} \times \cM_{\Sigma} \rightarrow
\cO,
\end{equation}
which is analogous to the one defined in \cite[1.5.3, 1.8.1]{dfg} 
in the case $F=\Q$ (note that since $\Sigma\supset P_\rho$ we do not
need the rather technical  \cite[Lemma 1.5]{dfg})

Theorem A is implied by the first part of the following :

\begin{theo} \label{main} 
Let $\rho: \G_F \rightarrow  \GL_2(\overline{\F}_p)$ be a
continuous representation satisfying
 $\mathrm{{\bf (LI}_{Ind \rho}{\bf )}}$ and
$\mathrm{{\bf (Mod}_{\rho}{\bf )}}$. Let $\Sigma$
be a finite set of primes containing $P_{\rho}$.  
Then $\pi_{\Sigma} : \cR_{\Sigma} \rightarrow
\cT_{\Sigma}$ is an isomorphism of complete intersections over $\cO$
and $\cM_{\Sigma}$ is free of finite rank over $\cT_{\Sigma}$.
In particular, all $\Sigma$-ramified deformations of $\rho$ are modular.

Moreover, for all Hilbert modular newforms $f$  such that
$\rho_{f,p}$ is a $\Sigma$-ramified deformations of
$\rho$  :
 \begin{equation}
\#\rH^1_{\Sigma}(F,\Ad^0(\rho_{f,p})\otimes
 \Q_p/\Z_p)=\#(\cO/\eta_{f}^{\Sigma})<\infty.
\end{equation}
\end{theo}

\begin{proof} We proceed by induction on $\#\Sigma$.
Assume first that $\Sigma=P_{\rho}$.  We already know that $\pi_{P_{\rho}} :
\cR_{P_{\rho} }\rightarrow \cT_{P_{\rho} }$ is an isomorphism of complete
intersections over $\cO$ and  
$\cM_P:=\cM$ is free of rank $2^d$ over 
$\cT_{P_{\rho} }$ (\cf Theorem \ref{minimal} if $P_{\rho}=\varnothing$
and  proposition \ref{isom3} together with Fujiwara \cite[Theorem 9.1]{Fu}  in general).

\medskip

Assume now that the theorem holds for some $\Sigma\supset P_{\rho}$, that is to
say  $\pi_{\Sigma}: \cR_{\Sigma}\rightarrow \cT_{\Sigma}$ is an isomorphism
of complete intersections over $\cO$ and that $\cM_{\Sigma}$ is free
over $\cT_{\Sigma}$. In particular, we have
$\#\Phi_{f}^{\Sigma}=\#(\cO/\eta_{f}^{\Sigma})$, where $f$ is 
a newform contributing to $\cM$. 

Let  $\gq$ be a prime outside $\Sigma$ not dividing $p$.
Put $\Sigma'=\Sigma\cup\{\gq\}$.

It follows directly from  Proposition \ref{tangent} and  Definition \ref{selmer}    that :
$$\#\Phi_{f}^{\Sigma'}\leq\#\Phi_{f}^{\Sigma}\cdot \#\rH^0(F_{\gq},
(\Ad^0(\rho_{f,p})\otimes\Q_p/\Z_p)(1)).$$

By Theorem \ref{numerical} and Proposition \ref{raising}, the theorem will hold
for $\Sigma'$ if we  construct a surjective homomorphism
$\cT_{\Sigma'}\twoheadrightarrow \cT_{\Sigma}$ 
compatible with the surjections $\theta_{f}^{\Sigma}$ 
and $\theta_{f}^{\Sigma'}$ and a
$\cT_{\Sigma'}$-linear injective homomorphism
$\xi:\cM_{\Sigma}\hookrightarrow \cM_{\Sigma'}$ inducing a surjection 
$\widehat{\xi}: \cM_{\Sigma'}\twoheadrightarrow \cM_{\Sigma}$ such that
$\widehat{\xi}\circ\xi(\cM_{\Sigma})=T\cdot\cM_{\Sigma}$ for
some $T\in \cT_{\Sigma}$ satisfying
\begin{equation}\label{cong-ideal}
\#(\cO/\theta_{f}^{\Sigma}(T))=\#\rH^0(F_{\gq},(\Ad^0(\rho_{f,p})\otimes\Q_p/\Z_p)(1)).
\end{equation}

This is done case by case, depending on the local behavior of $\rho$
 at $\gq$ (\cf Definition \ref{P}).  

The case  $\gq\in \Sigma_{\rho}\bs
 S_{\rho}$ is relatively straightforward, since adding such a prime does not
change $\cM_{\Sigma}$. We will distinguish two more cases :

\medskip

1) Assume that $\gq\notin\Sigma_{\rho}$. In this
case  $\rho_{\gq}\otimes\nu_{\gq}^{-1}$ is  unramified.

By Theorem \ref{ihara} and
 Proposition \ref{ihara-bis}, the homomorphism
$$\pr_3^*\pr_1^*+\pr_3^*\pr_2^*+\pr_4^*\pr_1^* : \cM_{\Sigma}^{\oplus 3}=
\rH^d(Y_{\Sigma},\bV_{\cO})[\psi ,\widetilde{\nu}\phi]_{\gm_{\rho}}^{\oplus 3}
\rightarrow\rH^d(Y_{\Sigma'},\bV_{\cO})[\psi
,\widetilde{\nu}\phi]_{\gm_{\rho}}$$ 
is injective with flat cokernel.

The  characteristic polynomial  of $U''_{\gq}$
acting on $\cM_{\Sigma}^{\oplus 3}$ is
$X(X^2-T_{\gq}X+S_{\gq}\Nm(\gq))$ and  $X=0$ is  simple  root modulo
$\varpi$ of this polynomial.
 Hence the localization of  the above injection at
 $(U''_{\gq},U_{\gu}-\alpha_{\gu})$  yields 
 another injection with flat cokernel : 
$$\xi:\cM_{\Sigma}\overset{\sim}{\rightarrow}
(\cM_{\Sigma}^{\oplus 3})_{U''_{\gq}} \hookrightarrow
\cM_{\Sigma'}.$$
 This gives a surjective homomorphism
$\cT_{\Sigma'}\twoheadrightarrow (\cT_{\Sigma}^3)_{U''_{\gq}}\cong
\cT_{\Sigma}$.   Computations performed by Wiles \cite[\S2]{wiles2}
and Fujiwara \cite{Fu} 
show that $\widehat{\xi}\circ\xi(\cM_{\Sigma})=T\cdot\cM_{\Sigma}$ with
$$ T=(\Nm(\gq)-1)(T_{\gq}^2-S_{\gq}(\Nm(\gq)+1)^2).$$
Then (\ref{cong-ideal}) follows by a straightforward computation.

\medskip
2)  Assume that $\gq\in S_{\rho}$. In this
case  $\dim(\rho_{\gq}\otimes\nu_{\gq}^{-1})^{I_v}=1$.

By Proposition \ref{ihara-bis} there is  an exact sequence whose last arrow has a  flat cokernel :
$$ 0 \rightarrow \rH^d(Y_{K_{\Sigma}\cdot K'''_{\gq}},\bV_{\cO})[\psi ,\widetilde{\nu}\phi]_{\gm_{\rho}}
\overset{(\pr_1^*,-\pr_2^*)}{\longrightarrow}
\rH^d(Y_{\Sigma},\bV_{\cO})[\psi ,\widetilde{\nu}\phi]_{\gm_{\rho}}^{\oplus 2}
\overset{\pr_3^*+\pr_4^*}{\longrightarrow}\rH^d(Y_{\Sigma'},\bV_{\cO})[\psi ,\widetilde{\nu}\phi]_{\gm_{\rho}},$$
where $K'''_{\gq}=\ker(K_1(\gq^{c_{\gq}-1})\overset{\det}{\longrightarrow}
\go_{\gq}^\times \overset{\widetilde{\nu}_{\gq}\phi}{\longrightarrow} \cO^\times)$.

The  characteristic polynomial  of $U''_{\gq}$
acting on $(\pr_3^*+\pr_4^*)(\cM_{\Sigma}^{\oplus 2})$ is
$X(X-U'_{\gq})$ and   $X=0$ is  simple  root  modulo $\varpi$ of this polynomial.
Hence the localization of the map $\pr_1^*+\pr_2^*$ at
$\gm_{\Sigma'}=(\gm_{\Sigma}, U''_{\gq})$  yields an injection with flat cokernel :
$$\xi:\cM_{\Sigma}\overset{\sim}{\rightarrow}
(\cM_{\Sigma}^{\oplus 2})_{U''_{\gq}} \hookrightarrow
\cM_{\Sigma'}.$$ 
This gives a surjective homomorphism
$\cT_{\Sigma'}\twoheadrightarrow (\cT_{\Sigma}^2)_{U''_{\gq}}\cong
\cT_{\Sigma}$.   Computations performed by Wiles \cite[\S2]{wiles2}
and Fujiwara \cite{Fu} 
show that $\widehat{\xi}\circ\xi(\cM_{\Sigma})=T\cdot\cM_{\Sigma}$ with
$$T=\begin{cases}
\Nm(\gq)-1 \text{ if } \rho_{\gq} \text{ is decomposable, and }
\\\Nm(\gq)^2-1 \text{ if } \rho_{\gq} \text{ is
  indecomposable}.\end{cases}$$ 
As above, (\ref{cong-ideal}) is obtained  by a straightforward computation.
\end{proof}

\subsection{Towards the modularity of a quintic threefold.}
 
We will give now an example coming from the geometry where Theorem A
applies.   Consani and  Scholten \cite{CoSc} consider the middle
degree cohomology of a quintic threefold $\widetilde{X}$ (a proper and
smooth $\Z[\frac{1}{30}]$-scheme with Hodge numbers 
 $h^{3,0}=h^{2,1}=1$, $h^{2,0}=h^{1,0}=0$ and $h^{1,1}=141$). They
show that the $\G_{\Q}$-representation
$\rH^3(\widetilde{X}_{\overline{\Q}},\Q_p)$ is induced from a two
dimensional representation $\widetilde{\rho}$ of $\G_{\Q(\sqrt{5})}$ and
conjecture the modularity of  $\widetilde{\rho}$. As explained in
\cite{dim-dieu}, Theorem \ref{main} implies the  following proposition

\begin{prop} {\rm (Dieulefait-D.)} Assume  $p\geq 7$ and that
  $\widetilde{\rho}$ is congruent modulo $p$ to the $p$-adic Galois
  representation attached to  a Hilbert modular form on $\Q(\sqrt{5})$ of
weight $(2,4)$ and some prime to $p$ level.
 Then $\widetilde{\rho}$ is  modular and,  in particular, the
 $L$-function  associated to 
  $\rH^3(\widetilde{X}_{\overline{\Q}},\Q_p)$ has  an analytic continuation
to the whole complex plane and satisfies a functional equation.
\end{prop}

\section{Cardinality of the adjoint Selmer Group.} \label{adjoint}

In this section we give a proof of Theorem B.  It is enough to
establish  (i), since then the finiteness of
$\rH^1_{\mathrm{f}}(F,\Ad^0(\rho_{f,p})\otimes \Q_p/\Z_p)$
implies (ii) by the same argument as  in \cite[\S2.2]{dfg}.

Choose a   finite  set $\Sigma$ of primes not dividing $p$, containing
the auxiliary prime $\gu$ and all primes (not dividing $p$) at which
$\Ad^0(\rho_{f,p})$ is 
ramified. Denote by $f_\Sigma$ the automorphic form contributing 
to $\cM_\Sigma$ corresponding to the newform $f$ of Theorem B.

\subsection{Periods of automorphic forms.}

For $J\subset J_F$ denote by $\epsilon_J$ the corresponding 
character  of the Weyl group $\left(\begin{smallmatrix} 1 &
    0 \\   0 & \pm 1 \end{smallmatrix}\right)^{J_F}\subset \GL_2(F\otimes_{\Q}\R)$. 
Put  $F_\infty=\left(\begin{smallmatrix} 1 &
    0 \\   0 & -1 \end{smallmatrix}\right)^{J_F}\in \GL_2(F\otimes_{\Q}\R)$.

We fix an isomorphism $\C\cong\overline{\Q}_p$
extending the embedding $\iota_p: \overline{\Q}\hookrightarrow
\overline{\Q}_p$.
Since $\cM_{\Sigma}$ is free over the 
principal domain $\cO$, lemma \ref{oldform}(ii) implies that 
$\cM_{\Sigma}[f,\epsilon_J]$
is a free $\cO$-module of rank one, where 
$[\epsilon_J]$ denotes the eigenspace for this character   and 
$[f]$ denotes  $\cap_{v\notin S}\ker(T_v-c(f,v))$.

Let $S_k(K_\Sigma;\psi ,\widetilde{\nu}\phi)$ be the $\C$-vector space of
automorphic forms on $\GL_2(F)\bs \GL_2(\A)$ which are  
holomorphic of weight $(k;k_0)$ at infinity  and right
$K_0(\gn_\Sigma)$-equivariant for the character 
$\left(\begin{smallmatrix} a & b \\ c &  d\end{smallmatrix}\right)\mapsto 
\psi (d)\widetilde{\nu}\phi(ad-bc)$. In 
particular such a form is right $K_\Sigma$-invariant. 
Let 
 \begin{equation} \begin{split}
\delta_J  & : S_k(K_\Sigma;\psi ,\widetilde{\nu}\phi)\rightarrow
 \rH^d_{\mathrm{cusp}}(Y_{\Sigma},\bV_{\C})[\psi
 ,\widetilde{\nu}\phi,\epsilon_J]   \text{ , and } \\ 
  \delta  & : \bigoplus_{J\subset J_F} S_{k,J}(K_\Sigma;\psi ,\widetilde{\nu}\phi)\rightarrow  
 \rH^d_{\mathrm{cusp}}(Y_{\Sigma},\bV_{\C})[\psi ,\widetilde{\nu}\phi].
\end{split}\end{equation}
 denote  the Matsushima-Shimura-Harder isomorphisms
(\cf \cite[Proposition 3.1, (4.2)]{hida74}). 

\begin{defin} \label{periods} For every $J\subset J_F$, we fix  a basis  $b_{f,J}$ of
$\cM_{\Sigma}[f,\epsilon_J]$ and  define the period
$\Omega_{f,J}=\frac{\delta_J(f_\Sigma)}{b_{f,J}}\in
\C^\times\!/\!\cO^\times$.  We fix  $J_0\!\subset\!  J_F$ and
put $\Omega_f^+=\Omega_{f,J_0}$ and  $\Omega_f^-= \Omega_{f,J_F\!\bs\!J_0}$.
\end{defin}

\begin{rque}
Classicaly, the Matsushima-Shimura-Harder periods of a newform $f$ of
level $\gn$ are defined using a basis of the free rank one $\cO$-module
$\rH^d(Y_1(\gn),\bV_{\cO})[f,\epsilon_J]$ (\cf \cite[\S4.2]{dim-hmv}).
As shown in  \cite[Theorem 6.6, \S4.4,\S4.5]{dim-hmv} the value at $1$ of the
imprimitive adjoint $L$-function divided by those periods 
measures the congruences modulo $p$ between $f$ and 
other Hilbert modular eigenforms of same  weight, level and
central character. However, in general the corresponding local 
Hecke algebra does not have a Galois theoretic interpretation 
whereas, as proved in Theorem A, $\cT_{\Sigma}$ does, hence our choice
to define  the periods using  $\cM_{\Sigma}[f,\epsilon_J]$. 
\end{rque}

Next we explain the relation between  the Peterson inner product 
and the modified Poincaré pairing defined in \S\ref{poincare} under the Matsushima-Shimura-Harder isomorphism. 

The Atkin-Lehner involution 
  $\iota=\left(\begin{smallmatrix} 0 &
    -1 \\ \gn_\Sigma &   0 \end{smallmatrix}\right)$ 
induces an isomorphism 
\begin{equation}
S_k(K_\Sigma;\psi ,\widetilde{\nu}\phi)
 \overset{\sim}{\longrightarrow} 
 S_k(K_\Sigma;\psi^{-1},(\widetilde{\nu}\phi)^{-1})
\quad f\mapsto f(\cdot \iota)\otimes  \psi^{-1}.
\end{equation}

The Hecke operator $[K_{\Sigma}xK_{\Sigma}]$ acts on 
$S_k(K_\Sigma;\psi ,\widetilde{\nu}\phi)$ by  sending  $f$ on 
$\sum_i f(\cdot x_i^{-1})$, where $K_{\Sigma}xK_{\Sigma}=\coprod_i
K_{\Sigma}x_i$.

We have the following commutative diagram, where the horizontal 
arrows are isomorphisms : 
\begin{equation}
\xymatrix@C=60pt@R=20pt{ 
S_k(K_\Sigma;\psi ,\widetilde{\nu}\phi) \ar[r]^{\iota^{-1}}  
\ar[d]^{[K_\Sigma \iota^{-1} x \iota K_\Sigma ]} & 
S_k(K_\Sigma^\iota;\psi ,\psi\phi^{-1}\widetilde{\nu}^{-1})
\ar[r]^{\otimes \psi^{-1}}
 \ar[d]^{[K_\Sigma^\iota x K_\Sigma^\iota]}&
S_k(K_\Sigma;\psi^{-1},(\widetilde{\nu}\phi)^{-1})
\ar[d]^{\psi (\det(x))[K_\Sigma xK_\Sigma ]} \\
S_k(K_\Sigma;\psi ,\widetilde{\nu}\phi) \ar[r]^{\iota^{-1}}  & 
S_k(K_\Sigma^\iota;\psi ,\psi\phi^{-1}\widetilde{\nu}^{-1})
\ar[r]^{\otimes \psi^{-1}} &
 S_k(K_\Sigma;\psi^{-1},(\widetilde{\nu}\phi)^{-1}) } 
\end{equation}

Finally, for $f_1,f_2\in S_k(K_\Sigma;\psi ,\widetilde{\nu}\phi)$ we define the
normalized Peterson inner  product by  
\begin{equation} \label{norm-peterson}
(f_1,f_2) =[K(1):K_0(\gn_\Sigma)]^{-1}
\int_{Y^{\ad}_\Sigma} f_1(g)\overline{f_2(g)}|\det(g)|^{2-k_0}dg.
\end{equation}  

We have $(f_1|_{[K_\Sigma xK_\Sigma]},f_2)=
|\det(x)|^{2-k_0}(f_1,f_2|_{[K_\Sigma x^{-1}K_\Sigma]})=
\psi (\det(x))(f_1,f_2|_{[K_\Sigma\iota x \iota^{-1}K_\Sigma]})$.

It follows that  the Hecke eigenvalues of $f_\Sigma$ are complex
conjugates of those of  $f_\Sigma(\cdot \iota)\otimes  \psi^{-1}$. 
Using  \ref{oldform}(ii) we deduce by   Strong Multiplicity One that 
these two forms differ by a constant, which 
turns out to be in $\cO^\times$ (the arguments of \cite[2.13]{dfg} 
involving local epsilon factors  can be adapted to our setting). 
Hence, in the computation that follows, this constant
 can be ignored, as well as $\N(\gn_\Sigma)$ and  powers of $2$ : 
$$(f_\Sigma,f_\Sigma) \cO=[\delta(f_\Sigma),
\delta(f_\Sigma(\cdot \iota F_\infty)\otimes  \psi^{-1})]\cO=
\langle \delta(f_\Sigma), \delta(f_\Sigma(\cdot F_\infty))  \rangle\cO =
\langle  \delta_{J_0}(f_\Sigma),\delta_{J_F\!\bs\!J_0}(f_\Sigma)
\rangle\cO.$$ 

From here and  definition \ref{periods}  we obtain the relation we
have been looking for :  
\begin{equation} \label{f0}
\langle b_{f,J_0},b_{f,J_F\!\bs\!J_0} \rangle \cO
=\frac{(f_\Sigma,f_\Sigma)}{\Omega_f^+\Omega_f^-}\cO.
\end{equation}

\subsection{The Rankin-Selberg method.}

The Rankin-Selberg method relating the Peterson inner product to 
the value at $1$ of the adjoint $L$-function has been carried out
by Shimura for Hilbert modular newforms $f$ of level $K_1(\gn)$. Since
the level  
structures $K_{\Sigma}$ that we consider are more general, the
resulting formula in our case slightly differs from Shimura's. 
While Shimura's formula relates the Peterson inner
product of $f$ with  the imprimitive adjoint $L$-function, in 
our setting the  Peterson inner product of $f_\Sigma$ will be
related to adjoint $L$-function outside $\Sigma$. 
We follow Jacquet's adelic version of the  Rankin-Selberg method
for $\GL_2$ and our main reference is Bump \cite{bump}.

All integrals that we consider are with respect to Haar measures on the corresponding
algebraic groups. The  normalized Peterson inner product
(\ref{norm-peterson}) can be rewritten as  : 

\begin{equation}
(f_\Sigma,f_\Sigma)=
\int_{\GL_2(F)\A^\times \bs \GL_2(\A)}|f_\Sigma(g)|^2|\det(g)|^{2-k_0}dg.\end{equation}

The automorphic form $f_\Sigma$ admits an adelic Fourier expansion : 
\begin{equation}
 f(g)=\sum_{y\in F^\times} W\left( \left(\begin{smallmatrix} y & 0 \\ 0 &
    1\end{smallmatrix}\right)g\right),\end{equation} 
where $W(g)=\int_{\A/F}\overline{\lambda}(x)
f\left(\left(\begin{smallmatrix} 1 & x \\ 0 & 1\end{smallmatrix}\right)g\right)dx$ is the
adelic Whittaker function with respect to an additive unitary character $\lambda$ given 
locally by $\lambda_v(x_v)=\exp(-2\pi i\tr(x))$. If $\delta_v$ denotes the valuation at $v$ of
the different $\gd$ of $F$, then $\lambda_v(v^{-\delta_v}\go_v)=1$. 
 The following decomposition can be found
in  \cite[Theorem 3.5.4]{bump}, but one should be careful to replace the usual $k/2$ by
$k+m-t$ since we are using the
arithmetic (non-unitary) normalization (\cf \cite[pp.566--567]{dimtildg}) :

\begin{equation}
W\left(\left(\begin{smallmatrix} y & 0 \\ 0 &
      1\end{smallmatrix}\right)\right)= 
y_{\infty}^{k+m-t}\exp(-2\pi\tr(y_{\infty}))\prod_v
W_v\left(\left(\begin{smallmatrix} y_v & 0 \\ 0 &
      1\end{smallmatrix}\right)\right).\end{equation}

Let  $\varphi$ be the
Schwartz function on $\A\times \A$ defined as product of the
following local functions : 
\begin{equation}
\varphi_\tau(x,y)=\exp(-\pi(x^2+y^2)) \text{ and }  
\varphi_v=\begin{cases}
\mathrm{char}(\go_v)\otimes\mathrm{char}(\go_v) & \text{ , for }v\notin \Sigma ;\\ 
\mathrm{char}(v^{c_v+d_v})\otimes\mathrm{char}(\go_v^\times)&\text{ ,
  for } v\in \Sigma.
\end{cases} \end{equation}

For $g\in\GL_2(\A)$ put $\varepsilon(g)=
\zeta_{F, \Sigma}(2s)^{-1}\pi^{sd}\Gamma(s)^{-d}
|\det(g)|^{s}\int_{\A_{F}^\times} 
|t|^{2s}\varphi(t(0,1)g)dt$.

Then $\varepsilon$ is a right $K_0(\gn_\Sigma)
\SO_2(F\otimes_{\Q}\R)$-invariant function  on $\GL_2(\A)$ such that $\varepsilon(1)=1$ and 
$\varepsilon\left(\left(\begin{smallmatrix} y & x \\ 0 & 
      y'\end{smallmatrix}\right)g\right)=
|\frac{y}{y'}|^{s}\varepsilon(g)$. Consider as in 
\cite[\S 3.7]{bump} the  Eisenstein series :
 \begin{equation}
E(g,s)=\sum_{B(F)\bs \GL_2(F)} \varepsilon(\gamma g). \end{equation}

The  Rankin-Selberg unfolding yields (cf \cite[pp.372--373]{bump}) : 
\begin{equation*}\begin{split}\int_{\GL_2(F)\A^\times \bs
      \GL_2(\A)} E(g,s) |f_\Sigma(g)|^2 |\det(g)|^{2-k_0}dg= \\
=\int_{B(\A) \bs \GL_2(\A)}\int_{\A_{F+}^\times}
\left|W\left(\left(\begin{smallmatrix} y & 0 \\ 0 &
        1\end{smallmatrix}\right)g\right)\right|^2 
\varepsilon\left(\left(\begin{smallmatrix} y & 0 \\ 0 &
      1\end{smallmatrix}\right)g\right) 
|y|^{1-k_0}|\det(g)|^{2-k_0}dydg.\end{split}\end{equation*}
 Here $\A_{F+}^\times$
denote the subgroup of idèles with totally positive infinite part. 
In \cite{bump} the integration is over $\A_{F}^\times$ but this makes
no difference, since $\A_{F}^\times=\A_{F+}^\times F^\times $ and 
the adelic Fourier expansion of $f(g)$ is supported only by 
totally positive elements. Using Iwasawa decomposition 
$$\GL_2(\A)=B(\A)\GL_2(\go\otimes \widehat{\Z}) \SO_2(F\otimes_{\Q}\R),$$
and the right $\SO_2(F\otimes_{\Q}\R)$-invariance of the integrand, we
further rewrite this integral as  
$\prod_\tau Z_\tau \prod_v Z_v$, 

$$ \text{ where  }Z_v =\int_{\GL_2(\go_v)}\int_{F_v^\times}
\left|W_v\left(\left(\begin{smallmatrix} y & 0 \\ 0 &
        1\end{smallmatrix}\right)g\right)\right|^2 
\varepsilon_v\left(\left(\begin{smallmatrix} y & 0 \\ 0 &
      1\end{smallmatrix}\right)g\right)
|y|^{1-k_0}dydg\text{ , and }$$
$$  Z_\tau=\int_0^{\infty} \left|W_\tau\left(\begin{smallmatrix} y & 0 \\ 0
      & 1\end{smallmatrix}\right)\right|^2 
\varepsilon_\tau\left(\begin{smallmatrix} y & 0 \\ 0 &
    1\end{smallmatrix}\right) |y|^{1-k_0}d^\times y =\int_0^{\infty}
\exp(-4\pi y) y^{s+k_\tau-1} d^\times y 
=(4\pi)^{-s-k_\tau+1}\Gamma(s+k_\tau-1).$$

Furthermore for $v\notin \Sigma$ (resp. $v\in \Sigma$) the function 
$W_v\cdot(\phi_v\widetilde{\nu}_v)^{-1}\circ \det$ is right $\GL_2(\go_v)$-invariant
(resp. $K_0(v^{c_v+d_v})$-invariant) and hence the same holds for $|W_v|^2$. 
Moreover $\varepsilon_v\left(\left(\begin{smallmatrix} y & 0 \\ 0 &
     1\end{smallmatrix}\right)g\right)|y|^{-s}$ is by definition 
the characteristic function of $\GL_2(\go_v)$
(resp. $K_0(v^{c_v+d_v})$). Hence for all $v$ :
$$Z_v=\int_{F_v^\times} \left|W_v\left(\begin{smallmatrix} y & 0 \\ 0 &
      1\end{smallmatrix}\right)\right|^2 |y|^{s+1-k_0} dy.$$

For $v\notin \Sigma$ we have 
$Z_v=\N(v^{\delta_v})^s(1+\N(v)^{-s}) L_v(\Ad^0(\rho_{f,p}),s)$ (\cf \cite[Proposition 3.8.1]{bump}).

For $v\in  \Sigma $, $W_v$ is  killed by $U_v$, hence
   $ Z_v=\N(v^{\delta_v})^s$. Therefore  
$$ \int_{\GL_2(F)\A^\times \bs \GL_2(\A)}E(g,s)|f_\Sigma(g)|^2
|\det(g)|^{2-k_0}dg= \frac{\N(\gd)^{s}\rL_{\Sigma }(\Ad^0(\rho_{f,p}),s)}
{\zeta_{F,\Sigma}(2s)\zeta_{F,\Sigma}(s)^{-1}}
\prod_{\tau\in J_F} \frac{\Gamma(s+k_\tau-1)}{(4\pi)^{s+k_\tau-1}}.$$

By \cite[Proposition 3.7.5]{bump}, $E(g,s)$
has a pole at $s=1$ with residue independent of $g$ and equal to 
the residue at $s=1$ of the function $\zeta_{F,\Sigma}(2)^{-1}\pi^d 
\int_{\A}\int_{\A^\times}|t|^{2s}\varphi(t,tx)dtdx$. One 
readily computes : 
$$\int_{\R}\int_{\R^\times}|t|^{2}\varphi(t,tx)dtdx=
\frac{1}{\pi}\int_{\R}\frac{dx}{(1+x^2)}=1
\text{ , and }$$
$$ \int_{F_v}\int_{F_v^\times}|t|^{2s}\varphi(t,tx)dtdx=\begin{cases}
 (1-\N(v)^{1-2s})^{-1} & \text{ , for } v\notin \Sigma  ;\\
(1-\N(v)^{-1})(1-\N(v)^{1-2s})^{-1}\N(v)^{(1-2s)(c_v+d_v)}  & \text{ , for } v\in \Sigma. 
\end{cases} $$
\begin{equation}(f_\Sigma,f_\Sigma)  = \frac{\N(\gn_\Sigma\gd)}{2^{|k|}}
\Gamma(\Ad^0(\rho_{f,p}),1) \rL_{\Sigma}(\Ad^0(\rho_{f,p}),1)=
\frac{\rL_{\Sigma}(\Ad^0(\rho_{f,p}),1)}{\pi^{|k|+d}}
 \frac{ \prod_\tau (k_\tau-1)!}{4^{|k|} \N(\gn_\Sigma\gd)^{-1}}.
\end{equation}
 
  Since by our assumptions 
 $\frac{ \prod_\tau (k_\tau-1)!}{4^{|k|}\N(\gn_\Sigma\gd)^{-1}} \in
 \overline{\Z}^\times_{(p)}$ it follows that 
 $\frac{\rL_{\Sigma}(\Ad^0(\rho_{f,p}),1)}{\pi^{|k|+d}(f_\Sigma,f_\Sigma) }  
\in \overline{\Z}^\times_{(p)}$.  Since by definition 
$\Gamma(\Ad^0(\rho_{f,p}),s)=\prod_{\tau\in   
  J_F}\pi^{-(s+1)/2}\Gamma\left(\frac{s+1}{2}\right)
(2\pi)^{-(s+k_{\tau}-1)}\Gamma(s+k_{\tau}-1)$ we obtain 
\begin{equation} \label{f1}
\frac{\Gamma(\Ad^0(\rho_{f,p}),1)\rL_{\Sigma}(\Ad^0(
\rho_{f,p}),1)}{(f_\Sigma,f_\Sigma)} \in \overline{\Z}^\times_{(p)}.  
\end{equation}

\subsection{End  of the proof of Theorem B(ii).}

Recall that  $\cM_{\Sigma}$
is endowed with a perfect  $\cT_{\Sigma}$-linear
pairing : $\langle\cdot,\cdot\rangle_{\Sigma}:\cM_{\Sigma} \times \cM_{\Sigma}
\rightarrow \cO$.

Since for all $J\!\subset\!  J_F$, $\cM_{\Sigma}[f,\epsilon_J]$ is free
of rank one over $\cT_{\Sigma}$ it follows that : 
$$(\eta_{f}^{\Sigma})^2=
\mathrm{disc}(\cM_{\Sigma}[f,\epsilon_{J_0}]\oplus
\cM_{\Sigma}[f,\epsilon_{J_F\!\bs\!J_0}])=
\langle b_{f,J_0},b_{f,J_F\!\bs\!J_0} \rangle ^2 \cO.$$ 

Using (\ref{f0}) we obtain 
\begin{equation} \label{f2}
\eta_{f}^{\Sigma}=\langle b_{f,J_0},b_{f,J_F\!\bs\!J_0} \rangle \cO
=\frac{(f_\Sigma,f_\Sigma)}{\Omega_f^+\Omega_f^-}\cO.
\end{equation}

By $\mathrm{{\bf (LI}_{Ind \rho}{\bf )}}$,    
$\rho|_{\G_{F(\zeta_p)}}$ is irreducible, hence Schur's lemma imply
$\rH^0(F,\Ad^0(\rho_{f,p})\otimes \Q_p/\Z_p)=
\rH^0(F,(\Ad^0(\rho_{f,p})\otimes \Q_p/\Z_p)(1))=0$.
Then  \cite[Lemma 2.1]{dfg}, which remains valid over $F$,  yields :
$$\Fitt_{\cO}\left(\rH^1_{\Sigma}(F,\Ad^0(\rho_{f,p})
\otimes \Q_p\!/\!\Z_p)\right)=
\Fitt_{\cO}\left(\rH^1_{\mathrm{f}}(F,\Ad^0(\rho_{f,p})
\otimes\Q_p\!/\!\Z_p)\right)   \prod_{v\in\Sigma}  
\Fitt_{\cO}\left(\rH^1_{\mathrm{f}}(F_v,\Ad^0(\rho_{f,p})(1))^*\right).$$

By  \cite[Proposition I.4.2.2(i)]{FoPR} and \cite[p.708, Lemma 2.16]{dfg}
$$\Tam(\Ad^0(\rho_{f,p}))=\prod_\tau \Tam_\tau(\Ad^0(\rho_{f,p}))
\prod_v
\Tam_v(\Ad^0(\rho_{f,p}))=\prod_{v\in\Sigma}\Tam_\tau(\Ad^0(\rho_{f,p})).$$ 
Furthermore, by \cite[Proposition I.4.2.2(ii)]{FoPR} and its proof and
\cite[(57)]{dfg}, for $v\in\Sigma$ we have 
$$ \Tam_v(\Ad^0(\rho_{f,p}))= 
\Fitt_{\cO}\left(\rH^1(I_v,\Ad^0(\rho_{f,p}))^{\G_{F_v}}_{\mathrm{tor}}\right)=
\Fitt_{\cO}\left(\rH^1(I_v,\Ad^0(\rho_{f,p})(1))^{
\G_{F_v}}_{\mathrm{tor}}\right)= $$
$$=\frac{\Fitt_{\cO}\left(\rH^1_{\mathrm{f}}
(F_v,\Ad^0(\rho_{f,p})(1))\right)}{
\Fitt_{\cO}\left(\rH^1_{\mathrm{f}}(F_v,\Ad^0(\rho_{f,p})(1)
\otimes\Q_p)\right)}=  \rL_v(\Ad^0(\rho_{f,p}),1)
\Fitt_{\cO}\left(\rH^1_{\mathrm{f}}(F_v,\Ad^0(\rho_{f,p})(1))^*\right).$$

From the three previous equations  we deduce : 
\begin{equation}\begin{split} \label{f3}
\Tam(\Ad^0(\rho_{f,p})) \Fitt_{\cO}\left(\rH^1_{\mathrm{f}}(F,\Ad^0(\rho_{f,p})
\otimes \Q_p\!/\!\Z_p)\right)=\\
=\prod_{v\in\Sigma}  \rL_v(\Ad^0(\rho_{f,p}),1)
\Fitt_{\cO}\left(\rH^1_{\Sigma}(F,\Ad^0(\rho_{f,p})\otimes\Q_p\!/\!\Z_p)\right).
\end{split}\end{equation}

Finally, since $\rho_{f,p}$
is a $\Sigma$-ramified deformation of $\rho=\overline{\rho}_{f,p}$ (\cf
Definition \ref{global}) and 
$\Sigma\supset P_{\rho}$ (\cf Definition \ref{P}),  Theorem \ref{main}
yields 
\begin{equation} \label{f4}
\Fitt_{\cO}\left(\rH^1_{\Sigma}(F,\Ad^0(\rho_{f,p})\otimes\Q_p/\Z_p\right)=
\eta_{f}^{\Sigma}.
\end{equation}
The theorem results by putting together the equations (\ref{f1}), (\ref{f2}), (\ref{f3})
and (\ref{f4}).

\begin{acknowledgements}
I started  this  project a few years ago as I was working at the
University Paris 13 and I would like to thank Jacques Tilouine for his
interest and his helpful comments which enriched my work. 
Most of the results were written down at the University Paris 7
and I would like thank my colleagues there, and most particularly Loïc Merel, 
for many useful conversations. Finally, I would like to thank 
very heartily the referee for the careful reading of the manuscript and his 
constructive criticism.  His insightful comments completely  transformed 
the earlier version of the paper. Needless to say, responsibility for
inaccuracies or errors is entirely my own. 
\end{acknowledgements}

\end{document}